\documentclass{amsart}

\usepackage{amsmath, amssymb, amsthm}
\usepackage{graphicx}
\usepackage[top=3cm, bottom=1.75cm, left=2.5cm, right=2.5cm]{geometry}
\usepackage[colorlinks=true, linkcolor=blue, citecolor = blue]{hyperref}
\usepackage{mathrsfs}
\usepackage{wasysym}
\usepackage{color}
\usepackage{stmaryrd}

\newtheorem{thm}{Theorem}[section]

\newtheorem{assumpt}{Assumption}

\newtheorem{rmk}{Remark}[section]

\newcommand{\ubar}[1]{\underline{#1}}

\newcommand{\eps}{\varepsilon}
\newcommand{\expt}{\mathbb{E}}
\newcommand{\mbbP}{\mathbb{P}}

\newcommand{\R}{\mathbb{R}}

\newcommand{\om}{\omega}
\newcommand{\C}{\mathcal{C}}

\newcommand{\Th}{\hat{T}}
\newcommand{\Ind}{\mathbf{1}_{\{0\}}}
\newcommand{\la}{\langle}
\newcommand{\ra}{\rangle}

\newcommand{\Psiz}{\hat{\Psi}}

\newcommand{\pj}{\Pi}
\newcommand{\pjeps}{\pj^\eps}
\newcommand{\Xeps}{X^\eps}
\newcommand{\icond}{\varphi}

\newcommand{\nstsp}{\mathbf{M}}
\newcommand{\noisegent}{\mathfrak{G}}
\newcommand{\gnoise}{\xi}
\newcommand{\gnoiseeps}{\gnoise^\eps}

\newcommand{\etavg}{\eta^{\hlev}}
\newcommand{\ham}{\mathfrak{h}}
\newcommand{\hlev}{\hbar}

\newcommand{\hprc}{\mathcal{H}}
\newcommand{\avgHproc}{\check{h}}
\newcommand{\hprcdrift}{b}
\newcommand{\hprcdiff}{\sigma}

\newcommand{\igen}{\mathcal{L}}

\newcommand{\mkK}{\mathcal{K}} 
\newcommand{\mfc}{\mathfrak{c}}
\newcommand{\mfs}{\mathfrak{s}}

\newcommand{\llb}{\llbracket}
\newcommand{\rrb}{\rrbracket}

\title[Perturbations of critical DDE]{Perturbations of linear delay differential equations at the verge of instability}
\date{19 February, 2016.}
\author{N. Lingala}
\email{nlingala1@gmail.com}
\author{N. Sri Namachchivaya}
\address{University of Illinois, Urbana, IL, USA.}

\begin{document}

\begin{abstract}
The characteristic equation for a linear delay differential equation (DDE) has countably infinite roots on the complex plane. This paper considers linear DDEs that are on the verge of instability, i.e. a pair of roots of the characteristic equation lie on the imaginary axis of the complex plane, and all other roots have negative real parts. It  is shown that, when small noise perturbations are present, the probability distribution of the dynamics can be approximated by the probability distribution of certain \emph{one dimensional} stochastic differential equation (SDE) \emph{without delay}. This is advantageous because equations without delay are easier to simulate and one-dimensional SDE are analytically tractable. When the perturbations are also linear, it  is shown that the stability depends on a specific complex number. The theory is applied to study oscillators with delayed feedback. Some errors in other articles that use multiscale approach are pointed out.
\end{abstract}

\maketitle


\section{Introduction}\label{sec:intro}

Delay differential equations (DDE) arise when the evolution of a variable at any time depends on the history of the variable. The evolution of many physical systems depends on their history owing to finite conduction velocities. Naturally, these systems are modeled by DDE. DDEs arise in many areas: biological systems, population dynamics, machining processes, viscoelasticity, laser optics etc. See \cite{KolMysh} for description of some examples. Many models of physiological systems, disease models, population dynamics involve DDE---see Mackey-Glass equation \cite{sch_MackGlass} for example.

The subject of this paper is linear DDE at the verge of instability. For example, consider the equation
\begin{align}\label{eq:PREpap:exampleverge}
\dot{x}(t)=\kappa x(t-1).
\end{align}
Seeking a solution of the form $x(t)=e^{t\lambda}$, we find that $\lambda$ must satisfy the characteristic equation $\lambda - \kappa e^{-\lambda}=0$. When $\kappa \in (-\frac{\pi}{2},0)$, all roots of the characteristic equation have negative real parts (see corollary 3.3 on page 53 of \cite{Stepanbook}). When $\kappa=-\frac{\pi}{2}$ a pair of roots $\pm i\frac{\pi}{2}$ are on the imaginary axis and all others have negative real parts. When $\kappa<-\frac{\pi}{2}$ some of the roots have positive real part. Hence, the system \eqref{eq:PREpap:exampleverge} is on the verge of instability at $\kappa=-\frac{\pi}{2}$. We study effect of perturbations on such systems, for example, 
$$\dot{x}(t)= \left(-\frac{\pi}{2} +  \eps \gnoise(t)\right)x(t-1) $$
where $\gnoise$ is a noise and $\eps \ll 1$ denoting the strength of the perturbation.

Such instability situations arise, for example, in machining processes. An oscillator of the form 
\begin{equation}\label{eq:PREpap:chatterex}
\ddot{q}(t)+2\zeta \dot{q}(t)+p^2q(t)=-\kappa p^2\left[q(t)-q(t-r)\right]
\end{equation}
is used to describe a phenomenon called  `regenerative chatter'  in machining processes \cite{Gabor1}. The model is as follows: A cutting tool is placed on a workpiece that is attached to a shaft rotating with time period $r$. The tool vibrates as it cuts the material from the workpiece.  Let $q(t)$ describe the position of a point on the machine tool. The force acting on the tool is proportional to the depth of the chip being cut and the depth is approximated as the difference  
between the present position ($q(t)$) of the tool and its position one revolution earlier ($q(t-r)$). The coefficient $\kappa$ is the force coefficient which depends, among other factors, on the width of cut. It is known that, for a fixed $r$, there exists a critical $\kappa_c$ such that the amplitude $q$ of the oscillator decreases exponentially  if $\kappa < \kappa_c$ and increases exponentially if $\kappa > \kappa_c$. When $\kappa=\kappa_c$ oscillations  of constant amplitude persist. This oscillatory behavior is called `chatter'. In machining, the goal is to have a large rate of cut. The greater the rate, the larger is $\kappa$, and chatter occurs when $\kappa$ is larger than a critical value resulting in poor surface finish. Researchers explored the possibility of achieving chatter suppression by varying structual parameters of the tool like damping and stiffness (see \cite{boringbar_magnet}, \cite{chatter_review}). Suppose there are small random perturbations in the natural frequency $p$ in \eqref{eq:PREpap:chatterex} such that $p=p_o(1+\eps \sigma(\gnoise(t)))$ where $\sigma$ is a mean-zero function of the noise $\gnoise$ and $\eps \ll 1$ is the strength of the perturbation, then on expanding in powers of $\eps$ and discarding terms of higher order, we have
\begin{eqnarray}\label{perturb-p}
\ddot{q}(t)+2\zeta \dot{q}(t)+p_0^2q(t) &=& -\kappa p_0^2 \left[q(t)-q(t-r)\right]
\nonumber \\ &&
+ \eps\sigma(\gnoise(t))\left[-2(1+\kappa)p_0 q(t) \right]
+\eps\sigma(\gnoise(t))\left[2\kappa p_0 q(t-r)\right],
\end{eqnarray}
which can be studied as a perturbation of \eqref{eq:PREpap:chatterex}. Also, small random perturbations in the properties of the material being cut could affect the tool dynamics---see \cite{buck_kusk_noise_machine_tool}.

Delay equations on the verge of instability arise also, for example, in the study of eye pupil \cite{Long90}, and act of human balancing \cite{Yao_delay_SDDE_stand}. In \cite{PhysRevE.85.056214}, authors make a case for studying effect of noise on oscillators with delayed feedback. As a prototypical oscillator they consider the van der Pol model
\begin{align}
\ddot{q}(t)+\om_0^2q(t)+\eta q(t-r)=&\beta \dot{q}(t)+\kappa \dot{q}(t-r)  -bq^2(t)\dot{q}(t)+q(t)\xi(t) \label{eq:Gaud_vdPosc}
\end{align}
with $\xi$ a Gaussian white noise with zero mean and variance $\la \xi(t)\xi(t')\ra=2D\delta(t-t')$.

Deterministic and stochastic DDE have been well studied in literature---see for example the books \cite{Halebook} (deterministic) and \cite{SEAM_book} (stochastic). Deterministic DDE \emph{at the verge of instability} are also well studied---see \cite{ChowParet} for averaging approach, \cite{DasChat_det_multiple_scale} and \cite{Nayfeh_Delay_Multiscale} for multiscale approach. Stochastic DDE \emph{at the verge of instability}, with noise being white, are studied by employing multiscale approach in \cite{Klosek_Kuske_multi_siam}, \cite{Kuske_SD} and \cite{PhysRevE.85.056214}, \cite{Gaud_hopf_DDE_multi_scale};  by averaging approach in \cite{Fofan02_sampstabchatter}, \cite{Fofana_PEM_chatter}, \cite{Caihong}; and by center-manifold approach in \cite{HuttEPLosc}.

However, \cite{Klosek_Kuske_multi_siam}, \cite{Kuske_SD}, \cite{PhysRevE.85.056214}, \cite{Gaud_hopf_DDE_multi_scale}
have committed serious errors in the analysis. These are pointed out in the appendix \ref{appsec:PREpap:errorsofKuskeGaudFof}. Sections \ref{appsec:PREpap:errorsofKuskeGaudFof__Kusk} (errors of \cite{Klosek_Kuske_multi_siam}, \cite{Kuske_SD}) and \ref{appsec:PREpap:errorsofKuskeGaudFof__Gaud} (errors of \cite{PhysRevE.85.056214}, \cite{Gaud_hopf_DDE_multi_scale}) can be read without further preparation. However, to understand \ref{appsec:PREpap:errorsofKuskeGaudFof__Fof} (shortcomings of \cite{Fofan02_sampstabchatter}, \cite{Fofana_PEM_chatter}, \cite{Caihong}) the mathematical background in the later two sections would be needed. \cite{HuttEPLosc} considers stability of scalar delay systems with additive white noise but commit an error in their analysis---which would be pointed out in section \ref{sec:Discussion}. \cite{lefebvre} considers a different kind of instability (one root of characteristic equation is zero and all other roots have negative real parts), which is reviewed in section \ref{sec:Discussion}.

This article deals with systems that can be studied as perturbations of linear DDE at the verge of instability.  In recent articles \cite{LNNSNarxSD} and \cite{LNarxNonLip} we have shown rigorously that, under certain conditions, the dynamics of such systems forced by white noise can be approximated (in a distributional sense) by the dynamics of a {\bf \emph{one-dimensional}} stochastic differential equation (SDE) {\bf \emph{without delay}}. The purpose of this article is three-fold: 
\begin{enumerate}
\item To exploit the results of  \cite{LNNSNarxSD} and \cite{LNarxNonLip} to show how the analysis of systems at the verge of instability can be simplified. The advantage arises because equations without delay are easier to simulate and one-dimensional SDE are analytically tractable. The articles \cite{LNNSNarxSD} and \cite{LNarxNonLip} deal rigorously with scalar systems forced by white noise. In this article we give (without proofs) explicit formulas for the approximating dynamics of vector-valued systems forced by white noise (equations of the form \eqref{eq:detDDE_pert} and \eqref{eq:quadnon_main_in_short_form}).

The approach taken in this article is similar to that in \cite{Fofan02_sampstabchatter}, \cite{Fofana_PEM_chatter}, \cite{Caihong},  in the sense that all use the spectral theory for DDE and averaging. However, \cite{Fofan02_sampstabchatter}, \cite{Fofana_PEM_chatter}, \cite{Caihong} consider specific applications of the equations of the form \eqref{eq:detDDE_pert} but do not consider the stronger perturbations as in equation \eqref{eq:quadnon_main_in_short_form}. \cite{HuttEPLosc} also uses spectral theory for DDE, and deals with stronger perturbations in the scalar case using a center-manifold approach. When dealing with equation \eqref{eq:quadnon_main_in_short_form}, the averaging approach that we take does not assume the existence of center-manifold (rigorous results about center-manifold for stochastic DDE are not known\footnote{However see \cite{SEAMstabmanif} for related results. One of the special cases of theorem 4.1 of \cite{SEAMstabmanif} is the following: In the case that zero is a fixed point of a stochastic DDE and the stochastic system linearized about zero does not have zero as a lyapunov exponent then local stable and unstable manifolds exist. These manifolds are the set of initial conditions which converge to or diverge from zero at an exponential rate. }).  
Further, the formulas \eqref{eqn:explictcentmanterms_form_1}--\eqref{eqn:explictcentmanterms_form_2} presented here, regarding the stronger perturbations $G_q$ in  \eqref{eq:quadnon_main_in_short_form}, are of independent interest. When applied in the deterministic DDE setting, they provide an alternate way to compute the effect of center-manifold terms on the amplitude of critical mode (more details are provided in section \ref{sec:strongerpert}).

\item To point out the errors in existing approaches that deal with white noise case.
\item To study systems forced by other general kind of noises (for example a continuous-time two-state markov chain). Theoretical results for this case (equations of the form \eqref{eq:detDDE_pert_gennoise}) dealt in section \ref{sec:OTHnoises} do not appear anywhere else. A sketch of the proof of the main result (theorem \ref{thm:gennoiseweakconvg}) is provided in appendix \ref{appsec:sketchproofGenNoise}.
\end{enumerate}

These claims would become more clear after the next two sections where the mathematical framework is explained. Also, in the case where the perturbations are also linear, a complex number is identified which \emph{alone} dictates the stability of the system.


\section{Mathematical setup of DDE}\label{sec:DDEframework}


\subsection{Notation}
\begin{enumerate}
\item $e^{\lambda \bullet}$ means a function whose evaluation at $\theta \in \R$ is $e^{\lambda \theta}$ 
\item * as superscript indicates transpose, 
\item $\bar{z}$ is complex conjugate of $z$,  
\item $\ubar{v}\in \R^n$ means $\ubar{v}$ is $n\times 1$ matrix with each entry in $\R$ and $\ubar{v}\in \R^{n*}$ means $\ubar{v}$ is $1\times n$ matrix with each entry in $\R$. The line underneath serves as a reminder that the quantity is multidimensional. Similar for $\mathbb{C}^n$ and $\mathbb{C}^{n*}$.
\end{enumerate}


\subsection{Equations considered in the article}
Let $x(t)$ be a $\R^n$-valued process governed by a DDE with maximum delay $r$. The evolution of $x$ at each time $t$ requires the history of the process in the time interval $[t-r,t]$. So, the state space can be taken as $\C:=C([-r,0];\R^n)$, the space\footnote{The space $\C$ is Banach space when equipped with sup norm: $||\eta||:=\sup_{\theta \in [-r,0]}|\eta(\theta)|$ for $\eta \in \C$. } of continuous functions on the interval $[-r,0]$ with values in $\R^n$. At each time $t$, denote the $[t-r,t]$ segment of $x$ as $\pj_t x$, i.e. $\pj_t x \in \C$ and
$$\pj_t x(\theta)=x(t+\theta), \quad \text{ for }\theta \in [-r,0].$$
Now, a linear DDE can be represented in the following form
\begin{align}\label{eq:detDDE}
\begin{cases}\dot{x}(t)=L_0(\pj_t x), \qquad t\geq 0, \\
\pj_0x =\icond \in \C,
\end{cases}
\end{align}
where $L_0:\C\to \R^n$ is a continuous linear mapping on $\C$ and $\icond$ is the initial history required.  For example, $\dot{x}(t)=-\frac{\pi}{2}x(t-1)$ can be represented using the linear operator given by $L_0(\eta)=-\frac{\pi}{2}\eta(-1)$ for $\eta \in \C$.

We assume there exists a bounded matrix-valued function $\mu:[-r,0]\to \R^{n\times n}$, continuous from the left on the interval $(-r,0)$ and normalized with $\mu(0)=0_{n\times n}$, such that 
\begin{align}\label{eq:L0matrixrep}
L_0\eta = \int_{[-r,0]}d\mu(\theta)\eta(\theta), \quad \forall \eta \in \C.
\end{align}
This is not a restriction: every continuous linear operator $L_0$ has such a representation. For example, $\dot{x}=-\frac{\pi}{2}x(t-1)$ can be represented with $\mu(\theta)=\begin{cases}\frac{\pi}{2} \qquad \theta=-r, \\ 0 \qquad \theta > -r.\end{cases}$

This article deals with perturbations of linear DDE, i.e. equations of the form
\begin{align}\label{eq:detDDE_pert}
\begin{cases}dx(t)=L_0(\pj_t x)dt+\eps^2 G(\pj_t x)dt+\eps F(\pj_t x)dW(t), \quad t\geq 0, \\
\pj_0x =\icond \in \C,
\end{cases}
\end{align}
where $F,G:\C \to \R^n$ are possibly nonlinear, $W$ is $\R$-valued Wiener process and $\eps \ll 1$ is a small number signifying perturbation. The following equations are also considered:
\begin{align}\label{eq:detDDE_pert_gennoise}
\begin{cases}dx(t)=L_0(\pj_t x)dt+\eps^2 G(\pj_t x)dt+\eps \sigma(\gnoise(t)) F(\pj_t x)dt, \quad t\geq 0, \\
\pj_0x =\icond \in \C,
\end{cases}
\end{align}
where $F,G:\C \to \R^n$ are possibly nonlinear, $\gnoise$ is a noise process (satisfying some assumptions) and $\sigma$ is a mean-zero function of the noise $\gnoise$. For example, one can have $\gnoise$ as a finite-state markov chain.

As an example, consider $\dot{\tilde{x}}=\kappa \tilde{x}(t-1)-\tilde{x}^3(t)$ where $\kappa$ has small perturbations about $-\frac{\pi}{2}$ according to $\kappa=-\frac{\pi}{2} +\eps \sigma(\gnoise(t)) +\eps^2$ where $\gnoise$ is a noise. Then $x(t)=\eps^{-1} \tilde{x}(t)$ can be put in the form  \eqref{eq:detDDE_pert_gennoise} with $L_0(\eta)=-\frac{\pi}{2}\eta(-1)$, $F(\eta)=\eta(-1)$ and $G(\eta)=-\eta^3(0)+\eta(-1)$.

The operator $L_0$ is asumed to be such that the unperturbed system \eqref{eq:detDDE} is on the verge of instability, i.e. $L_0$ satisfies the following assumption.
\begin{assumpt}\label{ass:assumptondetsys}
Define $$\Delta(\lambda)\,=\,\lambda I_{n\times n}-\int_{[-r,0]}d\mu(\theta)e^{\lambda \theta},$$
where $I$ is the identity matrix.
The characteristic equation 
\begin{align}\label{eq:chareq}
det(\Delta(\lambda))=0, \qquad \lambda \in \mathbb{C}
\end{align}
has a pair of purely imaginary solutions $\pm i\om_c$ and all other solutions\footnote{Typically there are countably infinite other roots.} have negative real parts.
\end{assumpt}

Since \eqref{eq:detDDE_pert} and \eqref{eq:detDDE_pert_gennoise} would be studied as perturbations of the linear DDE \eqref{eq:detDDE}, a brief overview of the unperturbed system \eqref{eq:detDDE} would be given now.


\subsection{The unperturbed system \eqref{eq:detDDE}}

The content in this section can be found in chapter 7 of \cite{Halebook} and chapter 4 of \cite{Diekmanbook}.
\subsubsection{Projection onto eigenspaces}
The space $\C$ can be split as $\C=P\oplus Q$ where $P$ is the eigenspace of the critical eigenvalues $\pm i\om_c$. Since $P$ corresponds to the critical eigenvalues $\pm i\om_c$, the projection of the dynamics of the unperturbed system onto $P$ is purely oscillatory with frequency $\om_c$. Since $Q$ corresponds to the eigenvalues with negative real part, the projection of the dynamics of the unperturbed system onto $Q$ decays exponentially fast.

Here we show, given an $\eta\in \C$, how to find the projection  onto the space $P$. For details, see chapter 7 of \cite{Halebook} and chapter 4 of \cite{Diekmanbook}.

  Any $\eta \in \C$ can be written as $\eta = \pi\eta + (I-\pi)\eta$ where $\pi \eta \,\in\,P $ and $(I-\pi)\eta \,\in \, Q$. Here $\pi$ is the projection operator $\pi:\C \to P$ and $I$ is the identity operator. The projection $\pi$ can be constructed as follows: Let 
\begin{align}\label{eq:chosingPhi}
\Phi=[\Phi_1,\,\,\,\Phi_2], \quad \Phi_1(\bullet)=\ubar{d}e^{i\om_c \bullet}, \quad \Phi_2(\bullet)=\bar{\ubar{d}}e^{-i\om_c \bullet}
\end{align}
where $\ubar{d}\in \mathbb{C}^n$ is chosen such that 
\begin{align}\label{eq:chosingPhi_const}
\Delta(i\om)\,\ubar{d}=0_{n\times 1}.
\end{align}
Note that each $\Phi_i$ belongs to $C([-r,0];\mathbb{C}^n)$. Define the bilinear form  $\la \cdot,\cdot\ra :C([0,r];\mathbb{C}^{n*})\times C([-r,0],\mathbb{C}^n) \to \mathbb{C}$, given by
\begin{equation}\label{eq:bilinform}
\la \psi,\eta \ra := \psi(0)\eta(0)-\int_{-r}^0\int_0^\theta \psi(s-\theta)d\mu(\theta)\eta(s)ds.
\end{equation}
Let 
\begin{align}\label{eq:chosingPsi}
\Psi=\left[\begin{array}{c}\Psi_1 \\ \Psi_2\end{array}\right], \quad \Psi_1(\bullet)=c\,\ubar{d_2}e^{-i\om_c \bullet}, \quad \Psi_2(\bullet)=\bar{c}\,\bar{\ubar{d_2}}e^{i\om_c \bullet},
\end{align}
where $\ubar{d_2}\in \mathbb{C}^{n*}$ is chosen such that 
\begin{align}\label{eq:chosingPsi_const}
\ubar{d_2}\,\Delta(i\om)=0_{1\times n}
\end{align}
and the constant $c$ is chosen such that 
\begin{align}\label{eq:PsiPhieqDelta}
\la \Psi_i,\Phi_j\ra =\delta_{ij}.
\end{align}
(Here $\delta_{ij}=1$ if $i=j$ and zero if $i\neq j$.) 

Writing $\la \Psi,\eta\ra = \left[\begin{array}{c}\la \Psi_1,\eta\ra \\ \la \Psi_2,\eta\ra \end{array}\right]$ 
we obtain for the projection $\pi: \C \to P$,
\begin{align}\label{eq:projoper_intermsof_bilform}
\pi(\eta)= \Phi \la \Psi,\eta \ra  = \Phi_1\la \Psi_1,\eta\ra + \Phi_2\la  \Psi_2, \eta \ra.
\end{align}
Note that $\la \Psi_1,\eta\ra $ and $\la \Psi_2,\eta\ra $ are complex conjugates and so are $\Phi_1$ and $\Phi_2$.


\subsubsection{Behaviour of solution on the eigenspaces}
The solution to the unperturbed system \eqref{eq:detDDE} can be written as $$\pj_t x=\pi\pj_t x+ (I-\pi)\pj_tx=\Phi z(t) + y_t$$ where $z(t)=\la \Psi, \pj_t x\ra$ and $y_t=\pj_t x-\Phi z(t)$. Note that $z\in \mathbb{C}^2$ is a 2-component vector with $z_2=\bar{z_1}$, and $\Phi z (t)\in P$ and $y_t\in Q$. It can be shown that
\begin{equation}\label{eq:proj_eqns_unpert}
\dot{z}(t)=Bz(t), \qquad B=\left[\begin{array}{cc}i\om_c & 0 \\ 0 & -i\om_c \end{array}\right],
\end{equation} 
i.e. $z$ oscillate with constant amplitude and frequency $\om_c$. So, $2z_1z_2$ is a constant in time. Further, it can be shown that $||y_t||$ decreases\footnote{This is the sup norm on $\C$.} to zero exponentially fast (because the dynamics on $Q$ is governed by eigenvalues with negative real parts).


\subsection{The perturbed systems \eqref{eq:detDDE_pert} and \eqref{eq:detDDE_pert_gennoise}}
Define the function $\ham:\C\to \R$ by
\begin{align}\label{eq:hamfuncdef}
\ham(\eta):=2\la \Psi_1,\eta \ra \la \Psi_2,\eta \ra, \qquad \eta \in \C.
\end{align}
As noted above,
$$2z_1(t)z_2(t)=2\la \Psi_1,\pj_t x\ra \la \Psi_2, \pj_t x \ra=\ham(\pj_t x)$$ is a constant for the \emph{unperturbed} system \eqref{eq:detDDE}. When we deal with the perturbed system \eqref{eq:detDDE_pert} or \eqref{eq:detDDE_pert_gennoise}, the quantity $\hprc(t):=\ham(\pj_tx)$ evolves much slowly compared to $x$ and $z_i$. In \eqref{eq:detDDE_pert}, because a Weiner process has the property that `\emph{the rescaled process $t\mapsto \eps W(t/\eps^2)$ has the same probability distribution as that of a Wiener process}', the noise perturbations take $O(1/\eps^2)$ time to significantly affect the $\hprc$ dynamics. Also, the prturbation $G$ is of strength $\eps^2$. Hence, significant changes in $\hprc$ occurs only in times of order $1/\eps^2$. In \eqref{eq:detDDE_pert_gennoise}, even though the strength of the noise perturbation is $\eps$, because $\sigma$ is a mean-zero function of the noise, significant changes in $\hprc$ occurs only in times of order $1/\eps^2$.

Our claim is that, under certain conditions on the coefficients $F$ and $G$, the probability distribution of the process $\hprc(t/\eps^2)$ converges to the probability distribution of a SDE \emph{without delay}. 
Because of the nature of decay on $Q$, $||y_t||$ decays to small values exponentially fast, and so studying $\hprc$ is enough to obtain a good approximation to the behaviour of $x$ in \eqref{eq:detDDE_pert} and \eqref{eq:detDDE_pert_gennoise}. How to obtain the SDE is shown in later sections.

\begin{rmk}\label{rmk:whyHprcUseful}
The reason why studying $\hprc$ would be useful is the following: for the moment assume the part of solution in the stable eigenspace $Q$ is zero, i.e. $\pj_tx=\Phi z(t)$ and $(I-\pi)\pj_tx=0$. Then, for the $j^{th}$ component of $x$ we have $x_j(t)=(\pj_tx(0))_j=(\ubar{d})_jz_1(t)+(\bar{\ubar{d}})_jz_2(t)$ where $\ubar{d}$ is choosen in \eqref{eq:chosingPhi}. Noting that $z_2=\bar{z_1}$ and that dynamics of $z_i$ is predominantly oscillatory with frequency $\om_c$, we find that the dynamics of $x_j$ is predominantly oscillatory with amplitude $2|(\ubar{d})_jz_1|$ or what is the same $\sqrt{4(\ubar{d})_j(\bar{\ubar{d}})_jz_1z_2}=|(\ubar{d})_j|\sqrt{4z_1z_2}=|(\ubar{d})_j|\sqrt{2\hprc}.$ Hence the magnitude of $\hprc$ indicates the amplitude of oscillation of $x$ (usually the amplitude might differ  from $|(\ubar{d})_j|\sqrt{2\hprc}$ by a slight amount because the part of the solution in $Q$, i.e. $(1-\pi)\pj_tx$ is not exactly zero).
\end{rmk}

A crucial role is played by the vector $\Psi(0)$. So the symbol $\Psiz$ is reserved for $\Psi(0)$.
\begin{align*}
\Psiz \overset{\mathrm{def}}{=\joinrel=} \Psi(0).
\end{align*}


\section{The perturbed system (\ref{eq:detDDE_pert})}

As noted above $\ham(\pj_t x)$ for the perturbed system \eqref{eq:detDDE_pert} varies slowly compared to $x$. Changes in $\ham(\pj_t x)$ are significant only on times of order $1/\eps^2$. Hence, we rescale time and write $X^\eps(t)=x(t/\eps^2)$ where $x$ is governed by \eqref{eq:detDDE_pert}. 

Under the above time-scaling, the $x$ time-series would be compressed by a factor of $\eps^2$. So, in order to be able to write the evolution equation for $X^\eps$, we need to define a new segment extractor $\pjeps_t$ as follows: for a $\R^n$ valued function $f$ defined on $[-\eps^2r,\infty)$ the $[t-\eps^2r,t]$ segment is given by
\begin{align}\label{eq:newsegextract}
(\pjeps_t f)\,(\theta)=f(t+\eps^2\theta), \qquad \quad -r\leq \theta \leq 0.
\end{align}

Now,  the process $X^\eps$ has the same probability law as that of a process satisfying
\begin{align}
dX^\eps(t)&=\frac{1}{\eps^2}L_0(\pjeps_t \Xeps)dt+ G(\pjeps_t \Xeps)dt  + F(\pjeps_t \Xeps)dW(t), \quad t\geq 0, \qquad \pjeps_0\Xeps =\icond \in \C, \label{eq:detDDE_pert_rescale}
\end{align}
where $W$ is $\R$-valued Wiener process\footnote{We have used the fact that for a Wiener process $W$, $\eps W(t/\eps^2)$ has the same probability law as a Wiener process.}.

Write $\hprc^{\eps}(t):=\ham(\pjeps_t \Xeps)$ with $\ham$ defined in \eqref{eq:hamfuncdef}. Using Ito formula, it can be shown that $\hprc^{\eps}(t)$ satisfies
\begin{align}\label{eq:evolofhprceps}
d\hprc^{\eps}(t)=\hprcdrift(\pjeps_t \Xeps)dt+\hprcdiff(\pjeps_t \Xeps)dW, \qquad \hprc^{\eps}(0)=\ham(\icond),
\end{align}
where 
\begin{align}
\hprcdrift(\eta)&=E(\eta) G(\eta)  +\,\frac12 4 (\Psiz_1 F(\eta))(\Psiz_2 F(\eta)), \label{eq:hprcdrifdef}\\
\hprcdiff(\eta)&=E(\eta) F(\eta), \label{eq:hprcdiffdef} \\
E(\eta)&=2(\la \Psi_1,\eta\ra \Psiz_2+\la \Psi_2,\eta\ra \Psiz_1). \label{eq:hprcdrifE} 
\end{align}

Recall that we can write the solution as $\pjeps_t\Xeps=\Phi z(t)+ (I-\pi)\pjeps_t\Xeps$ where $z(t):=\la \Psi,\pjeps_t\Xeps \ra$. Note that the evolution of $z_i(t)=\la \Psi_i,\pjeps_t\Xeps \ra$ is fast compared to the evolution of $\hprc^\eps$ and is predominantly oscillatory. Heuristically, the $z_i$ oscillate fast \emph{along} trajectories of constant $\ham$ (the effect of $\frac{1}{\eps^2}L_0$) while at the same time diffusing slowly \emph{across} the constant $\ham$ trajectories (the effect of perturbations $G, F$).  Hence, the $z_i$ in the above coefficients $\hprcdrift$ and $\hprcdiff$ can be averaged.

\begin{thm}\label{thm:PREpap:convgHforWnoquad}
In the case when \\
(i) $F$ is constant and $G$ has stabilizing effect or \\
(ii) $F$ is either linear or constant and $G$ is Lipschitz,\\
the probability distribution of $\hprc^{\eps}$ from \eqref{eq:evolofhprceps} until any finite time $T>0$ converges, as $\eps \to 0$, to the probability distribution of a process $\avgHproc$ which is the solution of the SDE
$$d\avgHproc(t)=b_H(\avgHproc(t))dt+\sigma_H(\avgHproc(t))dW(t), \qquad \quad \avgHproc(0)=\ham(\icond),$$
where $b_H$ and $\sigma_H$ are obtained by averaging the functions in \eqref{eq:hprcdrifdef} and \eqref{eq:hprcdiffdef} as described below in section \ref{subsec:PREpap:evalbHsigH}. The perturbation $G$ is said to have `stabilizing effect' if the deterministic system $\dot{\hlev}=b_{H}(\hlev)$ is stable.
\end{thm}

Note that $\hprc$ encodes information only about the critical component $\pi \pjeps \Xeps$ of the solution. The above results should be augmented with a result that the stable component $(I-\pi)\pjeps \Xeps$ is small. Proof of theorem \eqref{thm:PREpap:convgHforWnoquad} and a result to the effect that the stable component of the solution is small are  presented  in \cite{LNarxNonLip} (also see \cite{LNNSNarxSD} for the case when $G$ is Lipschitz and $F$ is constant).

\subsection{Evaluation of $b_H$ and $\sigma_H$}\label{subsec:PREpap:evalbHsigH}
To evaluate $b_H$ and $\sigma_H$ at a specific value $\hlev \in \R$, we consider a solution $\pj_t x$ of the \emph{unperturbed} system \eqref{eq:detDDE} that remains in the space $P$ for all time and such that $\ham(\pj_tx)=\hlev$. For this purpose define
\begin{align}\label{eq:unpertsolused4avg}
\etavg_t\,\,\overset{\text{def}}=\,\,\frac12\sqrt{2\hlev}\,\Phi \left[\begin{array}{c}e^{i\om_c t} \\ e^{-i\om_c t}\end{array}\right].
\end{align}
Note that $\etavg_t\in P$ for all time and the $z$ coordinates of $\etavg_t$ given by $\frac12\sqrt{2\hlev}\left[\begin{array}{c}e^{i\om_c t} \\ e^{-i\om_c t}\end{array}\right]$ evolve according \eqref{eq:proj_eqns_unpert}. Hence $\etavg_t$ is the solution of the unperturbed system with the initial condition $\etavg_0$. Further, $\ham(\etavg_t)=2(\frac12\sqrt{2\hlev}e^{i\om_ct})(\frac12\sqrt{2\hlev}e^{-i\om_ct})=\hlev$.

Now, the averaged coefficients $b_H$ and $\sigma_H$ are given by
\begin{align} \label{eq:avgbH}
b_{H}(\hlev)&=\frac{1}{2\pi/\om_c}\int_0^{2\pi/\om_c}b\left(\etavg_t\right)dt,\\ \label{eq:avgsigH}
\sigma_{H}^2(\hlev)&=\frac{1}{2\pi/\om_c}\int_0^{2\pi/\om_c}\sigma^2\left(\etavg_t\right)dt.
\end{align}
The following fact would be useful in the evaluation of above averages: for $\etavg_t$, $E$ defined in \eqref{eq:hprcdrifE} becomes (on using \eqref{eq:PsiPhieqDelta})
\begin{align*}
E(\etavg_t) = \sqrt{2\hlev}(\Psiz_1e^{-i\om_c t}+\Psiz_2e^{i\om_c t}).
\end{align*}


\section{Examples}

In this section we show three examples. The first is a simple scalar system---we study the perturbations of $\dot{x}(t)=-\frac{\pi}{2}x(t-1)$. In section \ref{sec:subsec:a_scalar_eq}, while studying cubic nonlinear perturbations and additive white noise perturbations, we illustrate the results of previous section and show how the averaged process can yield information about the $x$ process. This example is a running one in the sense that we revisit it when studying stronger deterministic perturbations in section \ref{sec:strongerpert} and different kinds of noise in section \ref{sec:OTHnoises}.

The purpose of the second example is to propose a conjecture. When perturbations are linear as well, we identify a complex number and claim that it \emph{alone} dictates the stability of the system. We provide support to our conjecture using numerical simulations on $\dot{x}(t)=-\frac{\pi}{2}x(t-1)$.

The third is the van der Pol oscillator \eqref{eq:Gaud_vdPosc}. Here we illustrate the stabilizing/destabilizing effects of noise and show how the averaging results obtained in the previous section give good enough description of the effects of noise and allow us to compute how much bifurcation thresholds are displaced in presence of noise when compared to the deterministic case.


\subsection{A scalar equation}\label{sec:subsec:a_scalar_eq}
Consider the following equation:
\begin{align}\label{eq:examplesys}
dx(t)=-\frac{\pi}{2}x(t-1)dt + \eps^2 x^3(t-1)dt + \eps \sigma dW.
\end{align}
In this case $L_0\eta=-\frac{\pi}{2}\eta(-1)$, $G(\eta)=\eta^3(-1)$ and $F(\eta)=\sigma$. The characteristic equation $\lambda+\frac{\pi}{2}e^{-\lambda}=0$ has countably infinite roots on the complex plane. The roots with the largest real part are $\pm i\om_c =\pm i\frac{\pi}{2}$. Let $\Phi(\theta)=[e^{i\frac{\pi}{2}\theta}\,\,\,e^{-i\frac{\pi}{2}\theta}]$. Now, $\Psi$ can be evaluated (using \eqref{eq:bilinform} to \eqref{eq:PsiPhieqDelta}) to be
$$\Psi(\bullet)=\left[\begin{array}{c}(1+i\frac{\pi}{2})^{-1}e^{-i\frac{\pi}{2} \bullet} \\ (1-i\frac{\pi}{2})^{-1}e^{i\frac{\pi}{2} \bullet} \end{array}\right].$$

The averaged drift and diffusions can be calculated using \eqref{eq:hprcdrifdef}--\eqref{eq:avgsigH}  as
\begin{align}
b_H(\hlev)&=2\Psiz_1\Psiz_2\sigma^2-\frac32(i(\Psiz_1-\Psiz_2))  \hlev^2, \label{eq:avgbHforscalarex}\\
\sigma_H^{2}(\hlev)&=4\Psiz_1\Psiz_2\sigma^2 \hlev. \label{eq:avgsigHforscalarex}
\end{align}

In section \ref{subsec:numverif}, we illustrate how the averaged equation $d\hlev=b_H(\hlev)dt+\sigma_H(\hlev)dW$ can be used to gain information about \eqref{eq:examplesys} (recall remark \ref{rmk:whyHprcUseful}). The section \ref{subsec:numverif} can be read now, setting $\gamma_q=0$ in \eqref{eq:examplesys_quadadded_forsimu}.


\subsection{Linear perturbations}\label{subsec:linearpertWhitenoise} 
In this section we consider the case where perturbations are also linear, and identify a complex number which \emph{alone} dictates the stability of the system. Note that we restrict to systems satisfying assumption \ref{ass:assumptondetsys}. \cite{SEAM_AOP_2} discusses methods to obtain bounds on the maximal exponential growth rates of more general class of delay equations. However the bounds given in \cite{SEAM_AOP_2} are not optimal for systems satisfying assumption \ref{ass:assumptondetsys}.

Consider
\begin{align}\label{mult:eq:examplesys}
dx(t)=L_0(\pj_t x)dt + \eps L_1(\pj_tx) dW(t),
\end{align}
where $L_i$ are linear operators, with $L_0$ satisfying assumption \ref{ass:assumptondetsys}. 
The averaged equation corresponding to \eqref{mult:eq:examplesys} is
\begin{align}\label{mult:eq:examplesysavgd}
d\hlev(t)\,=b_H(\hlev) dt\,+\,\sigma_H(\hlev)\,dW(t),
\end{align}
where $b_H$ and $\sigma_H$ can be evaluated using \eqref{eq:hprcdrifdef}--\eqref{eq:avgsigH}  as
\begin{align*}
b_H(\hlev)&=C_b\hlev, \qquad  \sigma_H^2(\hlev)=C_{\sigma}\hlev^2, \\
C_b&=(\Psiz_1L_1\Phi_1)(\Psiz_2L_1\Phi_2)+(\Psiz_1L_1\Phi_2)(\Psiz_2L_1\Phi_1), \\
C_{\sigma}&=(\Psiz_1L_1\Phi_1+\Psiz_2L_1\Phi_2)^2+2(\Psiz_1L_1\Phi_2)(\Psiz_2L_1\Phi_1).
\end{align*}
The solution to \eqref{mult:eq:examplesysavgd} is given by
\begin{align}\label{eq:PREpap:linpertsolavg}
\hlev(t)=\hlev(0)\exp\left((C_b-\frac12C_\sigma)t+\sqrt{C_\sigma}W(t)\right).
\end{align}
The Lyapunov exponent \emph{for the averaged equation} \eqref{mult:eq:examplesysavgd} can be calculated to be 
\begin{align}
\lambda_{avg}&=\lim_{t\to\infty}\frac{1}{t}\log \hlev(t)  \notag \\
&=\lim_{t\to\infty}\frac{1}{t}\log \hlev(0)+(C_b-\frac12C_\sigma)+\sqrt{C_{\sigma}}\lim_{t\to\infty}\frac{W(t)}{t} \notag \\
&=(C_b-\frac12C_\sigma) \notag \\
&=-\frac12\left((\Psiz_1L_1\Phi_1)^2+(\Psiz_2L_1\Phi_2)^2\right). \notag
\end{align}

Define $\lambda_j^\eps(t):=\frac1t \log\, \sup_{s\in[t-mr,t]}|x_j(s)|$ with $m\in \mathbb{N}$ such that $mr > \frac{2\pi}{\om_c}$ (here $m$ is chosen so as to avoid oscillations in the modulus of $x$). We conjecture that for large $t$, $\lambda^\eps(t)$ is close to $\eps^2\frac12\lambda_{avg}$. The $\frac12$ arises from the fact that $\hlev$ is quadratic in $x$.

\begin{figure}[h!]
\centering
  \includegraphics[scale=0.6]{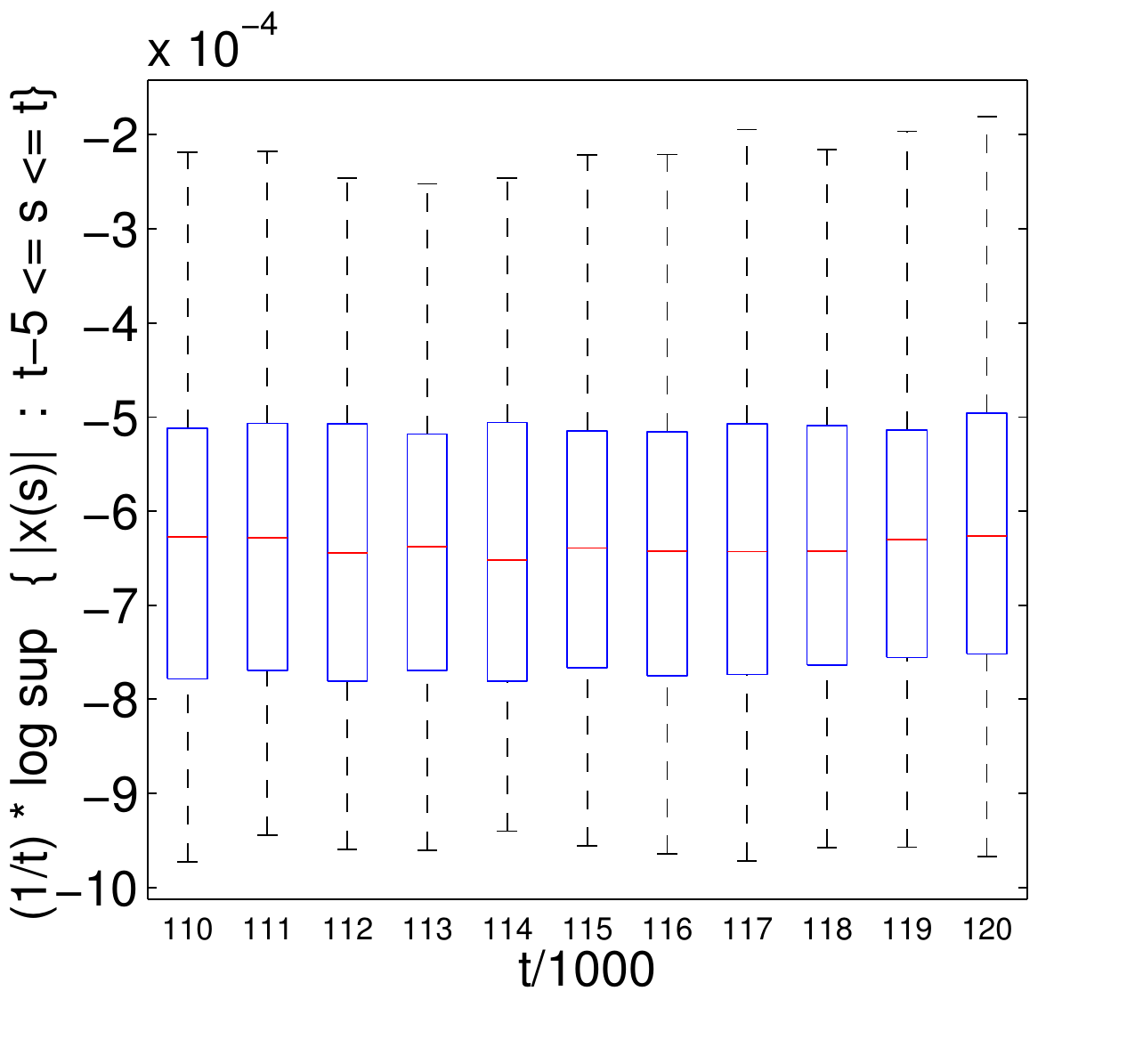}
  \caption{Box-plot of $\lambda^\eps(t):=\frac1t \log\, \sup_{s\in[t-5r,t]}|x(s)|$ for $t$ betwen 110,000 and 120,000 in steps of 1000. Red line is the mean of 80 realizations. Lower blue line is 25th percentile and upper blue line is 75th percentile.}
  \label{mult:fig:linearlyapfull_box}
\end{figure}

We verify the above conjecture using the sytem: 
\begin{align}\label{mult:eq:examplesys:particularex}
dx=-\frac{\pi}{2}x(t-1)dt+ \eps x(t-1)dW,
\end{align} 
i.e. $L_0\eta=-\frac{\pi}{2}\eta(-1)$ and $L_1\eta=\eta(-1)$. The Lyapunov exponent for \eqref{mult:eq:examplesysavgd} can be calculated to be $\lambda_{avg} \approx -0.122$ (the matrices $\Psiz$ and $\Phi$ are already calculated in section \ref{sec:subsec:a_scalar_eq}).
Eighty realizations of trajectories of \eqref{mult:eq:examplesys:particularex} are simulated with $\eps=0.1$ and initial condition $(\pj_0x)(\theta)=\cos(\om_c\theta)$ for $\theta \in [-r,0]$. In the figure 
\ref{mult:fig:linearlyapfull_box} we show the box plot
for $\lambda^\eps(t):=\frac1t \log\, \sup_{s\in[t-5,t]}|x(s)|$. For $t$ large, mean of $\lambda^\eps(t)$ is close to $-0.0006$ and we have $\eps^2\frac12\lambda_{avg} \approx -0.0006$. For details of the numerical scheme see appendix \ref{appsec:numericalscheme}.

Recalling that $\Psiz_2$ and $L_1\Phi_2$ are the complex conjugates of $\Psiz_1$ and $L_1\Phi_1$ respectively, we find that 
$$\lambda_{avg}=-\,Re[(\Psiz_1 L_1\Phi_1)^2]=-|\Psiz_1 L_1\Phi_1|^2\cos(2\theta_*),$$
where $\theta_*$ is the angle of the complex number $\Psiz_1L_1\Phi_1$. The stability condition $\lambda_{avg}<0$ translates to $\cos(2\theta_*)>0$. If the conjecture that for large $t$, $\lambda^\eps(t)$ is close to $\eps^2\frac12\lambda_{avg}$ is true, then the complex number $\Psiz_1L_1\Phi_1$ alone dictates the stability of \eqref{mult:eq:examplesys}.


\subsection{van der Pol oscillator}\label{subsec:vdpolosc}
In this section we consider the oscillator modeled by equation \eqref{eq:Gaud_vdPosc}, which was considered in \cite{PhysRevE.85.056214}. In studying \eqref{eq:Gaud_vdPosc}, our intentions are three fold: (i) to point out\footnote{This is done in appendix \ref{appsec:PREpap:errorsofKuskeGaudFof}} the errors in the analysis of \cite{PhysRevE.85.056214}, (ii) illustrate the stabilizing/destabilizing effects of noise, (iii) show that the averaging results obtained in the previous section give good enough description of the effects of noise.

The oscillator \eqref{eq:Gaud_vdPosc} has natural frequency $\om_0$ which would be altered by the delayed-feedbacks $\eta q(t-r)$ and $\kappa \dot{q}(t-r)$. Negative of $\beta$ indicates the strength of linear damping in the oscillator. The coefficient $b$, if positive, is the strength of nonlinear damping in the oscillator.

Since we intend to study the effect of small noise perturbations, we scale $D=\eps^2\tilde{D}$ with $\eps \ll 1$. Since we study the dynamics close to the zero fixed point, we zoom-in and write $x_1(t)=\frac{1}{\eps} q(t)$ and $x_2(t)=\frac{1}{\eps} \dot{q}(t)$. Then, the oscillator \eqref{eq:Gaud_vdPosc} can be put in the following form (using Ito interpretation)
\begin{align}\label{eq:Gaud_vdPosc_stp}
dx(t)=L_0(\pj_tx)dt &+\eps^2\left(\begin{array}{c}0 \\ -{b}x_1^2(t)x_2(t)\end{array}\right)dt  +\eps \sqrt{2\tilde{D}}\left(\begin{array}{c}0 \\ x_1(t)\end{array}\right)dW(t)
\end{align}
where $W$ is Wiener process and $L_0\phi=\int_{-r}^0d\mu(\theta)\phi(\theta)$ with
\begin{align*}
d\mu(\theta)=\left(\begin{array}{cc}0 & 1 \\ -\om_0^2 & \beta \end{array}\right)\delta_0(\theta)+\left(\begin{array}{cc}0 & 0 \\ -\eta & \kappa \end{array}\right)\delta_{-r}(\theta),
\end{align*}
where $\delta_0$ and $\delta_{-r}$ are delta functions, i.e. $\int \delta_0\phi = \phi(0)$ and $\int \delta_{-r}\phi = \phi(-r)$ for $\phi \in \C$.

\begin{figure}[h!]
\centering
\includegraphics[scale=0.7]{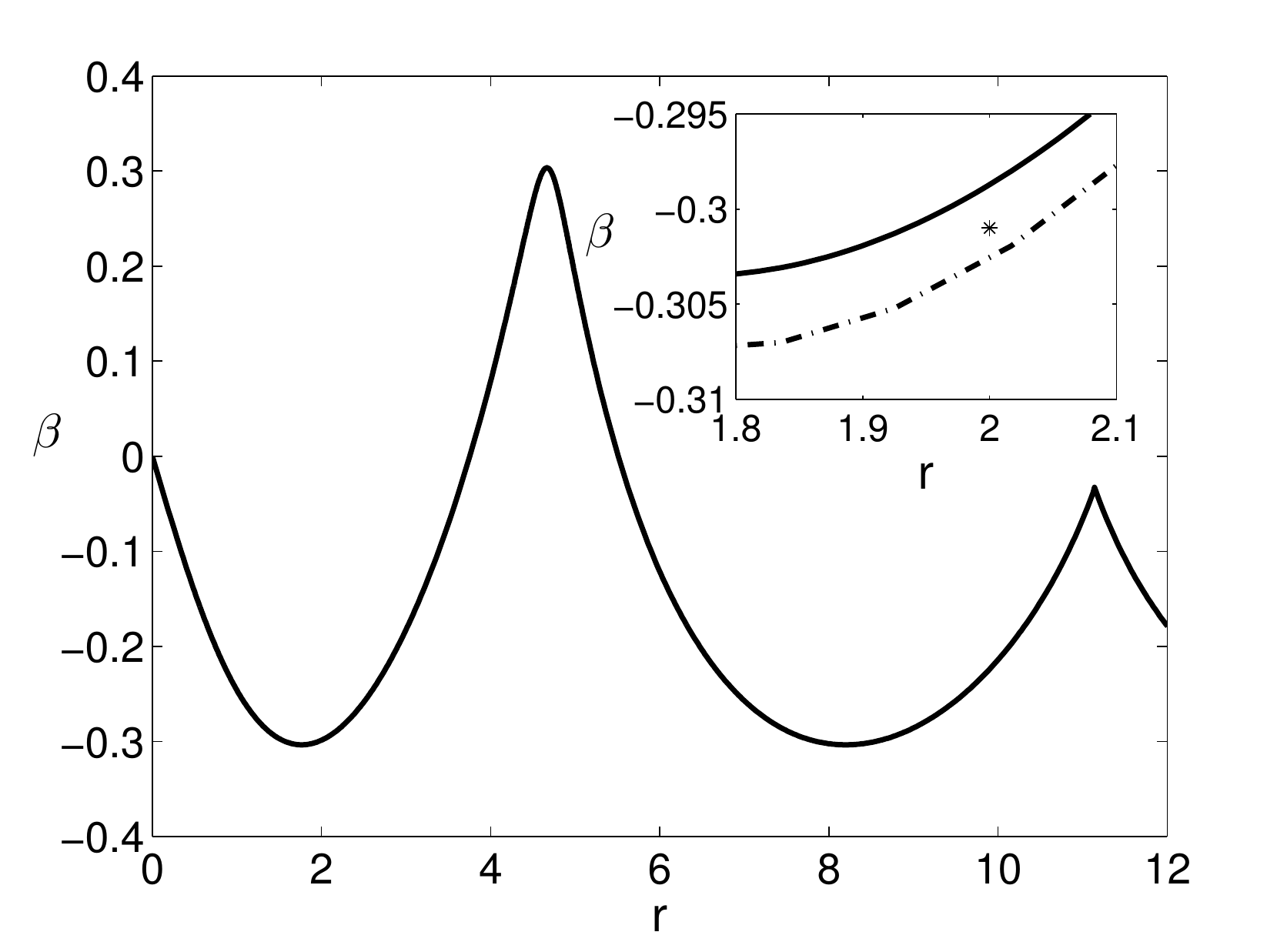}
  \caption{Boundary of stability for the fixed point $(x_1=0,x_2=0)$ of the system \eqref{eq:Gaud_vdPosc_stp} with $\eps=0$, $\om_0=1$, $\kappa=0$, $\eta=0.3$. For each delay $r$ there exists a critical value $\beta_c$ such that for $\beta<\beta_c$ the fixed point is stable and for $\beta>\beta_c$ the fixed point is unstable. In the inset, (theoretically predicted) stability boundary in presence of noise is shown with dashed line (obtained using \eqref{eq:criticalbetaforvanderpol_alt}). For this, $\eps=0.1$, $\tilde{D}=1$ and $b=1$. For $\beta$ in the region below the dashed line, theoretical results predict that  the $(0,0)$ fixed point is stable in presence of noise. Above the dashed line the fixed point looses stability; nevertheless invariant density exists. So, theoretical results predict that the noise has destabilized the region between solid and dashed lines. The point marked by $*$ in the inset is $r=2$, $\beta=-0.301$. For this point we show in figure \ref{fig:Gaud_InvDens_with_Noise} the invariant density obtained by numerical simulations. The theoretically obtained invariant density (obtained in \eqref{eq:GaudInvDens}) is in very good agreement with the actual density obtained from numerical simulations.}
  \label{fig:Gaud_StabBoundary}
\end{figure}

The characteristic equation becomes 
\begin{equation}\label{eq:vanderpolcharact}
-\lambda \beta +\lambda^2 +(\eta-\kappa\lambda) e^{-\lambda r}+\om_0^2=0.
\end{equation}
Since our intention is to study the effect of small noise perturbations on the oscillator when it is at the verge of instability, we assume that the parameters of the problem are such that the characteristic equation has two roots $\pm i\om_c$ on the imaginary axis and all other roots have negative real parts. With this assumption the unperturbed system $\dot{x}(t)=L_0(\pj_tx)$ is on the verge of instability. Figure \ref{fig:Gaud_StabBoundary} shows the stability boundary. 

 The matrices $\Phi$ and $\Psi$ can be evaluated (using \eqref{eq:chosingPhi} to \eqref{eq:PsiPhieqDelta}) as
\begin{align*}
\Phi(\bullet)&=\left(\begin{array}{cc}e^{i\om_c\bullet} & e^{-i\om_c\bullet} \\ i\om_c e^{i\om_c\bullet} & -i\om_c e^{-i\om_c\bullet}\end{array}\right)=\left(\begin{array}{c}\Phi_1 \,\,\,\Phi_2\end{array}\right), \\
\Psi(\bullet)&=\left(\begin{array}{cc}c(\om_0^2+\eta e^{-i\om_c r})e^{-i\om_c\bullet} & c(-i\om_c)e^{-i\om_c\bullet} \\ \bar{c}(\om_0^2+\eta e^{i\om_c r}) e^{i\om_c\bullet} & \bar{c}(i\om_c) e^{i\om_c\bullet}\end{array}\right)=\left(\begin{array}{c}\Psi_1\\ \Psi_2\end{array}\right),
\end{align*}
where 
\begin{align}\label{eq:constcforgaudcase}
c=(\om_c^2+e^{-i\om_c r}(\eta+i\eta r\om_c+\kappa r \om_c^2)+\om_0^2)^{-1}.
\end{align}
\begin{rmk}\label{rmk:hphyssignif}
The process $\ham(\pj_tx)$ with $\ham$ defined in \eqref{eq:hamfuncdef} has additional significance for this problem. If $\pj_tx$ was such that the stable part $(I-\pi)\pj_tx$ was zero, then $\pj_tx=\pi\pj_tx=\Phi z(t)$, which gives $$x(t)=\pj_tx(0)=\Phi_1(0)z_1(t)+\Phi_2(0)z_2(t)=\left[\begin{array}{c}z_1(t)+z_2(t)\\i\om_c(z_1(t)-z_2(t))\end{array}\right]$$
from which we get $\ham(\pj_tx)\overset{\text{by def}}=2z_1(t)z_2(t)=\frac12((x_1(t))^2+(x_2(t)/\om_c)^2)$ which represents some kind of energy in the oscillator (note that $x_1$ is position and $x_2$ is velocity). Usually $||(I-\pi)\pj_tx||$ decays to very small values exponentially fast and hence $\ham(\pj_tx)$ differs from the `energy' $\frac12((x_1(t))^2+(x_2(t)/\om_c)^2)$ by a little amount.
\end{rmk}
Using \eqref{eq:hprcdrifdef}--\eqref{eq:avgsigH} we have 
\begin{align*}
b_H(\hlev)&=(2\tilde{D})2|c|^2\om_c^2 \hlev -{b}\om_c^2\frac12(c+\bar{c})\hlev^2, \\
\sigma_H^2(\hlev)&=(2\tilde{D})\left(2|c|^2\om_c^2+(i\om_c(\bar{c}-c))^2\right)\hlev^2.
\end{align*}

To understand whether noise has a stabilizing or destabilizing effect, lets consider the damping $\beta$ as a bifurcation parameter. Write $\beta=\beta_c+\eps^2\tilde{\beta}$ and assume that at $\eps=0$, $\beta$ satisfies the characteristic equation \eqref{eq:vanderpolcharact}. Then, the effect of $\,\,\tilde{\beta}$ is to add another term $\tilde{\beta}(c+\bar{c})\om_c^2\hlev$ to $b_H$. Then, we can write the averaged equation as
\begin{align}\label{eq:vanderpolavgd}
d\hlev&=b_H(\hlev)dt+\sigma_H(\hlev)dW,
\end{align}
where
\begin{align*}
b_H(\hlev)&=C_b\hlev+C_b^{(2)}\hlev^2, \qquad \sigma_H^2(\hlev)=C_\sigma \hlev^2,
\end{align*}
$$C_b=(2\tilde{D})2|c|^2\om_c^2\left(1+\frac{\tilde{\beta}}{2\tilde{D}}\frac{(c+\bar{c})/2}{|c|^2}\right),$$
$$C_b^{(2)}=-b\om_c^2\frac12(c+\bar{c}),$$
$$C_\sigma=(2\tilde{D})2|c|^2\om_c^2\left(1+\frac{2((\bar{c}-c)/2i)^2}{|c|^2}\right).$$
To focus on the effect of noise, for the moment we ignore the nonlinearities by setting $b=0$ in \eqref{eq:Gaud_vdPosc_stp}. Corresponding averaged system then becomes
\begin{align}\label{eq:vanderpolavgd_linearized}
d\hlev=C_b\hlev + \sqrt{C_\sigma}\hlev dW.
\end{align}
The above system is unstable when\footnote{note that the solution is similar to \eqref{eq:PREpap:linpertsolavg}.} $C_b-\frac12C_\sigma>0$, i.e. when
\begin{align}\label{eq:criticalbetaforvanderpol}
\frac{\tilde{\beta}}{2\tilde{D}|c|}\frac{(c+\bar{c})/2}{|c|}>\frac{((\bar{c}-c)/2i)^2}{|c|^2}-\frac{1}{2}.
\end{align}

Let $\varsigma_1=\frac{(c+\bar{c})/2}{|c|}$ and $\varsigma_2=\left(\frac{((\bar{c}-c)/2i)^2}{|c|^2}-\frac{1}{2}\right).$ It can be shown\footnote{Note that $sign(\varsigma_1)=sign(\frac{c+\bar{c}}{c\bar{c}})=sign(\frac{1}{c}+\frac{1}{\bar{c}})$. Using \eqref{eq:constcforgaudcase} we have $$c^{-1}+(\bar{c})^{-1}=2(\om_c^2+\om_0^2)+\eta(e^{i\om_c r}+e^{-i\om_c r})+ir\om_c e^{-i\om_c r}(\eta-i\om_c k)-ir\om_c e^{i\om_c r}(\eta+i\om_c k).$$ Employing $\lambda=\pm i\om_c$ in the characteristic equation \eqref{eq:vanderpolcharact} we get, $$ir\om_c e^{-i\om_c r}(\eta-i\om_c k)-ir\om_c e^{i\om_c r}(\eta+i\om_c k)=-2\beta_c r\om_c^2,$$ $$2\eta(e^{i\om_c r}+e^{-i\om_c r})=(\om_c^2-\om_0^2)(e^{i\om_c r}+e^{-i\om_c r})^2+\beta_c i\om_c(e^{2i\om_c r}-e^{-2i\om_c r}).$$ Hence $c^{-1}+(\bar{c})^{-1}=2(\om_c^2+\om_0^2)+\frac12(\om_c^2-\om_0^2)(e^{i\om_c r}+e^{-i\om_c r})^2+\frac12\beta_c i\om_c(e^{2i\om_c r}-e^{-2i\om_c r})-2\beta_c r\om_c^2$ which can be simplified as $c^{-1}+(\bar{c})^{-1}=2\om_c^2(1+\cos^2\om_c r)+2\om_0^2(1-\cos^2\om_c r)-\beta_c\om_c(2r\om_c+\sin2\om_c r)$ which is positive if $\beta_c<0$.} that if $\beta_c<0$, then $\varsigma_1>0$.

Assume $\beta_c<0$. Then, \eqref{eq:criticalbetaforvanderpol} holds when
\begin{align}\label{eq:criticalbetaforvanderpol_alt}
\frac{\tilde{\beta}}{2\tilde{D}|c|}>\frac{\varsigma_2}{\varsigma_1}.
\end{align}

If noise was not present, i.e. $\tilde{D}=0$ in \eqref{eq:Gaud_vdPosc_stp}, then the $(x_1=0,x_2=0)$ fixed point of \eqref{eq:vanderpolavgd_linearized} would have been unstable for any $\tilde{\beta}>0$ (this is because $-\tilde{\beta}$ specifies how much additional damping is present in the system). If noise is present and $\varsigma_2>0$, then the $(x_1=0,x_2=0)$ fixed point of \eqref{eq:vanderpolavgd_linearized} is stable even for $0<\tilde{\beta}<{2\tilde{D}|c|\varsigma_2}/{\varsigma_1}$. So, noise has a stabilizing effect if $\varsigma_2>0$. 

Similar reasoning shows that the noise has destabilizing effect if $\varsigma_2<0$. If the noise was not present, then the $(x_1=0,x_2=0)$ fixed point of \eqref{eq:vanderpolavgd_linearized} would have been stable for any $\tilde{\beta}<0$. If noise is present and $\varsigma_2<0$, then \eqref{eq:vanderpolavgd_linearized} is unstable even for ${2\tilde{D}|c|\varsigma_2}/{\varsigma_1}<\tilde{\beta}<0$. So, noise has a destabilizing effect if $\varsigma_2<0$. This is the scenario presented in the inset of figure \ref{fig:Gaud_StabBoundary}.

The  stability of  \eqref{eq:Gaud_vdPosc_stp} when $b\neq 0$ depends on the stability of averaged nonlinear system \eqref{eq:vanderpolavgd}. However the theorem \ref{thm:PREpap:convgHforWnoquad} deals with only weak convergence of probability distributions and hence is not adequate to transfer the stability properties from the averaged system to the original system \eqref{eq:Gaud_vdPosc_stp}. Neverthelss we give an account of the stability of the averaged system \eqref{eq:vanderpolavgd}. When the nonlinearity is destabilizing, i.e. $C_b^{(2)}>0$, the system \eqref{eq:Gaud_vdPosc_stp} cannot be stable. When $C_b^{(2)}<0$ and $C_b-\frac12C_{\sigma}<0$ then the trivial solution $\hlev=0$ is the only equilibrium point of \eqref{eq:vanderpolavgd} and is stable. When $C_b^{(2)}<0$ and $C_b-\frac12C_{\sigma}>0$ the trivial solution  of \eqref{eq:vanderpolavgd} becomes unstable; nevertheless an invariant density exists. It is given by (obtained by solving steady-sate Fokker-Planck equation)
\begin{equation}\label{eq:GaudInvDens}
p(\hlev)=\frac{\chi^{\frac{2C_b}{C_{\sigma}}-1}}{\Gamma(\frac{2C_b}{C_{\sigma}}-1)}\,\hlev^{2(\frac{C_b}{C_{\sigma}}-1)}\,e^{- \hlev \chi}, \qquad \chi=2(-C_b^{(2)})/C_{\sigma},
\end{equation}
where $\Gamma$ is the Gamma function. 

\begin{figure}[h!]
\centering
\includegraphics[scale=0.6]{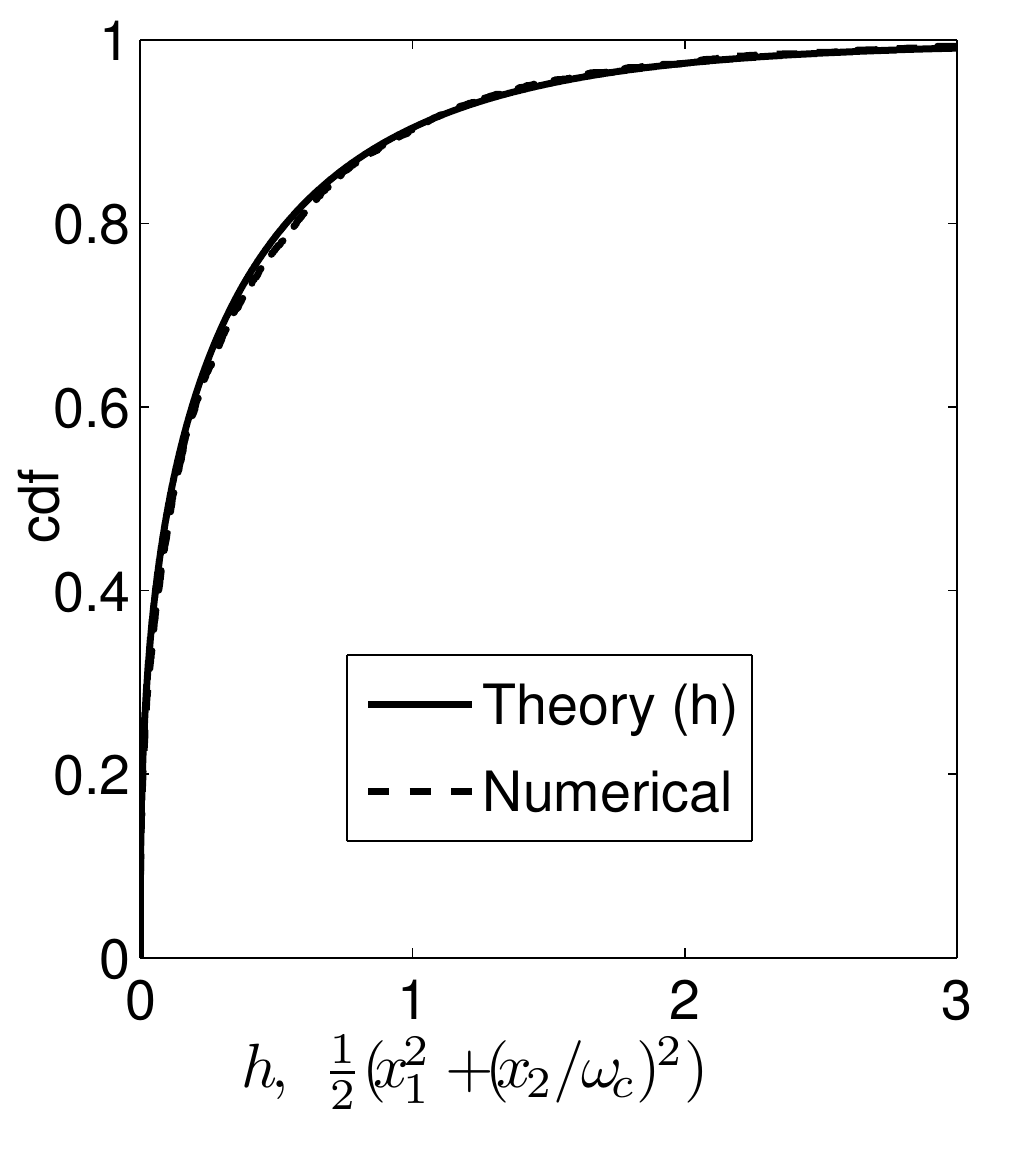}
  \caption{Cumulative distribution function (cdf) of the invariant density of $\frac12(x_1^2+(x_2/\om)^2)$ obtained from numerical simulation of \eqref{eq:Gaud_vdPosc_stp} with parameters specified by the point marked by `$*$' in the inset of figure \ref{fig:Gaud_StabBoundary} ($\om_0=1$, $\kappa=0$, $\eta=0.3$, $\eps=0.1$, $\tilde{D}=1$, $b=1$, $r=2$, $\beta=-0.301$). This agrees with the cdf of the density given in \eqref{eq:GaudInvDens}. For this case, the deterministic bifurcation threshold is $\beta_c=-0.2987$ and the predicted threshold in presence of noise is $\beta_c+\eps^2{2\tilde{D}|c|\varsigma_2}/{\varsigma_1}=-0.3027$.}
  \label{fig:Gaud_InvDens_with_Noise}
\end{figure}

The usefulness of the above results is shown in figure \ref{fig:Gaud_InvDens_with_Noise}. 
Let the parameters be specified by the point marked by `$*$' in the inset of figure \ref{fig:Gaud_StabBoundary}. 
When $\eps=0$, the $(x_1=0,x_2=0)$ fixed point of the oscillator \eqref{eq:Gaud_vdPosc_stp} would be stable because `$*$' lies below the stability boundary (solid line in figure \ref{fig:Gaud_StabBoundary}). However, in presence of noise the stability boundary is shifted by $\eps^2{2\tilde{D}|c|\varsigma_2}/{\varsigma_1}$ (dashed line in figure \ref{fig:Gaud_StabBoundary}). Now the fixed point loses stability; nevertheless invariant density exists. Numerical simulation is done with 3200 samples and the cumulative distribution function (cdf) of the invariant density of $\frac12(x_1^2+(x_2/\om_c)^2)$ is plotted in figure \ref{fig:Gaud_InvDens_with_Noise}. Also shown is the cdf arising from the averaging result \eqref{eq:GaudInvDens}. By the averaging theorems and remark \ref{rmk:hphyssignif} these two should be in good agreement---the figure \ref{fig:Gaud_InvDens_with_Noise} indeed shows this.

Numerical simulations in the case $\varsigma_2<0$ with $\eps=0.1$ show very good agreement with theoretical averaging results for $\beta$ in the range $\beta_c>\beta>\beta_c+0.9\eps^2({2\tilde{D}|c|\varsigma_2}/{\varsigma_1})$. Very close to the theoretically predicted bifurcation threshold in the presence of noise, i.e. $\beta \approx \beta_c+\eps^2({2\tilde{D}|c|\varsigma_2}/{\varsigma_1})$, the agreement is not very good. Actual bifurcation threshold in presence of noise (denoted by $\beta_{c,\text{noi}}$) obtained from numerical simulations of \eqref{eq:Gaud_vdPosc_stp}, is within $20\%$ of the theoretically predicted value, i.e. $\beta_c+\eps^2({2\tilde{D}|c|\varsigma_2}/{\varsigma_1})>\beta_{c,\text{noi}}>\beta_c+1.2\eps^2({2\tilde{D}|c|\varsigma_2}/{\varsigma_1})$. For details of the numerical scheme see appendix \ref{appsec:numericalscheme}.


\section{Stronger deterministic perturbations}\label{sec:strongerpert}
Here we consider systems with slightly stronger deterministic perturbations:
\begin{align}
dx(t)=L_0(\pj_t x)dt\,&+\,\eps G_q(\pj_t x)dt\,+\,\eps^2 G(\pj_t x)dt   +\,\eps F(\pj_t x)dW(t),  \label{eq:quadnon_main_in_short_form}
\end{align}
where $W$ is $\R$-valued Wiener process.

As an example, consider the noisy perturbation $d{\tilde{x}}=-\frac{\pi}{2} \tilde{x}(t-1)dt + \tilde{x}^2(t)dt + \eps^2\sigma dW$ of the DDE $\dot{\tilde{x}}(t)=-\frac{\pi}{2} \tilde{x}(t-1)+\tilde{x}^2(t)$. Then $x(t)=\eps^{-1} \tilde{x}(t)$ can be put in the form  \eqref{eq:quadnon_main_in_short_form} with $L_0(\eta)=-\frac{\pi}{2}\eta(-1)$, $F(\eta)=\sigma$, $G(\eta)=0$ and $G_q(\eta)=\eta^2(0)$.

The effect of $G_q$ in \eqref{eq:quadnon_main_in_short_form} is significant in just times of order $1/\eps$ whereas the effects of $G$ and $F$ are significant in times of order $1/\eps^2$. So we consider only those $G_q$ which are such that a certain kind of time averaged effect of $G_q$ is zero: 
\begin{align}\label{eq:assumption_on_Gq_zero}
\frac{1}{2\pi/\om}\int_0^{2\pi/\om}e^{-i\om_c t}\Psiz_1G_q(\etavg_t)\,dt=0,
\end{align}
where $\etavg_t$ is defined in \eqref{eq:unpertsolused4avg}. The assumption \ref{eq:assumption_on_Gq_zero} is a natural one: for example, $G_q$ which are homogenously quadratic in $\eta$ (say $G_q(\eta)=(\eta(0))^2$) satisfy the property \eqref{eq:assumption_on_Gq_zero}.

Writing $X^\eps(t)=x(t/\eps^2)$,  equation analogous to \eqref{eq:detDDE_pert_rescale} becomes
\begin{align}
dX^\eps(t)&=\frac{1}{\eps^2}L_0(\pjeps_t \Xeps)dt+ \frac{1}{\eps}G_q(\pjeps_t \Xeps)dt +G(\pjeps_t \Xeps)dt  + F(\pjeps_t \Xeps)dW(t), \quad t\geq 0, \label{eq:detDDE_pert_rescale_quadadded} \\
\pjeps_0\Xeps &=\icond \in \C. \notag
\end{align}
 Using Ito formula, $\hprc^{\eps}(t):=\ham(\pjeps_t \Xeps)$ satisfies
\begin{align}\label{eq:evolofhprceps_quad}
d\hprc^{\eps}(t)=\frac{1}{\eps}(\hprcdrift^{q,(1)}(\pjeps_t\Xeps)+ \hprcdrift^{q,(2)}(\pjeps_t\Xeps))dt+\hprcdrift(\pjeps_t \Xeps)dt+\hprcdiff(\pjeps_t \Xeps)dW, \qquad \hprc^{\eps}(0)=\ham(\icond),
\end{align}
where $\hprcdrift$, $\hprcdiff$ and $E$ are same as in \eqref{eq:hprcdrifdef}, \eqref{eq:hprcdiffdef}, \eqref{eq:hprcdrifE} respectively, and
\begin{align}
\hprcdrift^{q,(1)}(\eta)&=E(\eta)G_q(\pi\eta),  \label{eq:hprcdrifdef_quad1}\\ 
\hprcdrift^{q,(2)}(\eta)&=E(\eta)(G_q(\eta)-G_q(\pi\eta)). \label{eq:hprcdrifdef_quad2} 
\end{align}

Recall that we can write the solution as $\pjeps_t\Xeps=\Phi z(t)+ (I-\pi)\pjeps_t\Xeps$ where $z(t):=\la \Psi,\pjeps_t\Xeps \ra$. Note that the evolution of $z_i(t)=\la \Psi_i,\pjeps_t\Xeps \ra$ is fast compared to the evolution of $\hprc^\eps$ and is predominantly oscillatory. Heuristically, the $z_i$ oscillate fast \emph{along} trajectories of constant $\ham$ (the effect of $\frac{1}{\eps^2}L_0$) while at the same time diffusing slowly \emph{across} the constant $\ham$ trajectories (the effect of perturbations $G, G_q, F$).  Hence, the effect of $z_i$ in the above coefficients $\hprcdrift$ and $\hprcdiff$ can be averaged out. Our goal is to obtain an averaging result akin to theorem \ref{thm:PREpap:convgHforWnoquad}.
However, the terms arising from $G_q$ should be dealt with carefully. 
The assumption \ref{eq:assumption_on_Gq_zero} would entail that $\frac{1}{2\pi/\om}\int_0^{2\pi/\om}E(\etavg_t)G_q(\etavg_t)\,dt$ equals zero as well\footnote{This follows from the fact that $E(\etavg_t)=\sqrt{2\hlev}(e^{-i\om_c t}\Psiz_1+e^{i\om_c t}\Psiz_2)$ and $\Psiz_2$ is the conjugate of $\Psiz_1$.}.  Hence, when the oscillations are averaged, the leading order contribution of $\hprcdrift^{q,(1)}$ is zero. However, because of the $\frac{1}{\eps}$ multiplying $\hprcdrift^{q,(1)}$, higher order effects must be taken into account. 

We give explicit formulae for the contributions from $\hprcdrift^{q,(1)}$ and $\hprcdrift^{q,(2)}$, using solutions of the unperturbed system with $n$ specific initial conditions.  Atleast when $G_q$ is purely quadratic, \emph{the averaged terms arising from $\hprcdrift^{q,(k)}$ would be the same as what one gets from a formal center-manifold and normal-form calculation. However we do not assume the existence of a center-manifold.} The following method however has an advantage in that numerical integration can be used to find the answers. To provide an illustration of how the method works, a simple example without delay is worked in appendix \ref{appsec:strongdetpert_example_nodel}. To state the formulae, we need to set up some notation.


\subsection{Notation}\label{subsec:strongdetpert_not} 

 For $\icond \in \C$, let $\Th(t)\icond$ denote the solution at time $t$ of the unperturbed linear system \eqref{eq:detDDE} with initial condition $\pj_0x=\icond$, i.e. $\Th(t)\icond=\pj_t x$ where $x$ is governed by \eqref{eq:detDDE}.

Let $\Ind:[-r,0]\to \R^{n\times n}$ denote the matrix valued function 
\begin{align}\label{eq:matrixInd0}
\Ind(\theta)=\begin{cases} I_{n\times n}, \quad \theta =0, \\ 0_{n\times n}, \quad \theta \neq 0,\end{cases}
\end{align}
where $I$ is the identity matrix.
For a constant $n\times 1$ vector $\ubar{v}$, one can solve the unperturbed linear system \eqref{eq:detDDE} with $\pj_0x=\Ind \ubar{v}$. The solution is indicated by $\Th(t)\Ind\ubar{v}$.

 Recall that $\pi$ is the projection operator onto the critical eigenspace and is given by \eqref{eq:projoper_intermsof_bilform}. Even though $\Ind\ubar{v}$ does not belong to $\C$ (because it is not continuous), the definition $\pi(\Ind \ubar{v}):=\Phi \la \Psi, \Ind \ubar{v}\ra$ still makes sense\footnote{Rigorous way to extend the space $\C$ to include the discontinuities and the decomposition of the extended space as $P\oplus \hat{Q}$ is discussed in \cite{Halebook}.} using the bilinear form \eqref{eq:bilinform}.  On evaluation of the bilinear form we find that
\begin{align}\label{eq:aux:projofIndontoP}
\pi(\Ind \ubar{v})=\Phi \Psiz \ubar{v}.
\end{align}
 The meaning of $\Th(t)\pi\Ind \ubar{v}$ and $\Th(t)(I-\pi)\Ind \ubar{v}$ should now be clear. 

Suppose $G:\C \to \R^k$ and let $\eta,\xi \in \C$. Then $(\xi.\nabla)G(\eta)$ denotes the Frechet differential of $G$ evaluated at $\eta$ in the direction of $\xi$, i.e.
$$(\xi.\nabla)G(\eta)=\lim_{\delta \to 0}\frac{G(\eta+\delta \xi)-G(\eta)}{\delta}.$$
In a moment we would see the motivation for defining the following:
\begin{align}\label{eq:taudef}
\rho({\eta}):=\inf \left\{ t > 0\,:\,\la  \Psi, \Th(t)\pi\eta\ra = \frac12\sqrt{2\ham(\eta)}\left[\begin{array}{c}1\\1\end{array}\right] \right\},
\end{align}
\begin{align}\label{eq:a1qdef_aux}
a_q^{(1)}(\eta)=\int_0^{\rho(\eta)}\left(\left(\Th(s) \pi\Ind G_q(\eta)\right).\nabla\right)\hprcdrift^{q,(1)}(\Th(s)\pi \eta)ds,
\end{align}
\begin{align}\label{eq:a2qdef_aux}
a_q^{(2)}(\eta)&=\int_0^{\infty}\left(\left(\Th(s) \Ind G_q(\eta)\right).\nabla\right)\hprcdrift^{q,(2)}(\Th(s)\pi \eta)ds.
\end{align}


\subsection{Averaging}
\begin{thm}\label{thm:PREpap:convgHforWquad}
In the case when $F$ is constant and $G,G_q$ are Lipschitz; the probability distribution of $\hprc^{\eps}$ until any finite time $T>0$, converges as $\eps \to 0$, to the probability distribution of a process $\avgHproc$ which is the solution of the SDE
$$d\avgHproc(t)=(b_H+b_H^{q,(1)}+b_H^{q,(2)})(\avgHproc(t))dt +\sigma_H(\avgHproc(t))dW(t), \qquad \quad \avgHproc(0)=\ham(\icond),$$
where $b_H$ and $\sigma_H$ are same as in \eqref{eq:avgbH} and \eqref{eq:avgsigH} and $b_H^{q,(k)}$ for $k=1,2$ are given by 
\begin{align}\label{eq:avgbhqk}
b_H^{q,(k)}(\hlev)=\frac{1}{2\pi/\om_c}\int_0^{2\pi/\om_c}a_q^{(k)}\left(\etavg_t\right)dt,
\end{align}
where $\etavg_t$ is defined in \eqref{eq:unpertsolused4avg}.
\end{thm}
The proof of the above result can be found in \cite{LNNSNarxSD}. The key idea in obtaining the averaged effect of $G_q$ is this: Let $c^{q,(1)}$ be the function whose differential along the trajectory of the unperturbed system equals $\hprcdrift^{q,(1)}$ defined in \eqref{eq:hprcdrifdef_quad1}. Then the average effect of $\hprcdrift^{q,(1)}$ is negative of the average of `the differential of $c^{q,(1)}$ along the direction of the perturbations'. In symbols: the function $c^{q,(1)}(\eta)=-\int_0^{\rho(\eta)}\hprcdrift^{q,(1)}(\Th(s)\eta)ds$ is such that $\frac{d}{dt}\big|_{t=0}c^{q,(1)}(\Th(t)\eta)=\hprcdrift^{q,(1)}(\eta)$. The differential of $c^{q,(1)}$ along the direction of the perturbations is $(\Ind G_q(\eta).\nabla)c^{q,(1)}(\eta)$ which evaluates to $-a^{q,(1)}(\eta)$ (plus an additional term whose average turns out to be zero due to assumption \ref{eq:assumption_on_Gq_zero}). The average effect of $\hprcdrift^{q,(1)}$ is the average of $a^{q,(1)}$. Similar is the reasoning for $\hprcdrift^{q,(2)}$. For details see\footnote{\cite{LNNSNarxSD} deals with scalar systems and does not employ polar coordinates. Hence the form of expressions differ from here. However they evaluate to same numbers as here. The key difference is: \cite{LNNSNarxSD} writes an element $\eta \in P$ as $z_1\cos(\om_c \cdot)+z_2\sin(\om_c \cdot)$ with $z_i\in \R$. Here we write as $z_1e^{i\om_c\cdot}+z_2e^{-i\om_c \cdot}$ with $z_i\in \mathbb{C}$ and $z_2=\bar{z_1}$. } section 9 of \cite{LNNSNarxSD}. 
To illustrate the above idea, a simple example without delay is worked out in appendix \ref{appsec:strongdetpert_example_nodel}. We urge the reader to study appendix \ref{appsec:strongdetpert_example_nodel} to gain intuition about the process of obtaining the drift coefficients $b_H^{q,(i)}$.

The term $b_H^{q,(1)}$ is solely due to the critical eigenspace, and the term $b_H^{q,(2)}$ arises from the interaction between stable eigenspace and critical eigenspace. When $G_q$ is purely quadratic, these are the same terms that arise from a formal center-manifold calculation.

Note that $\hprc$ encodes information only about the critical component of the solution $\pi \pjeps \Xeps$. The above results should be augmented with a result that the stable component $(I-\pi)\pjeps \Xeps$ is small. Proof of theorem \ref{thm:PREpap:convgHforWquad} and a result to the effect that the stable component of the solution is small are presented in \cite{LNNSNarxSD}.

\begin{rmk}\label{rmk:Foferror1}
It is clear from \eqref{eq:hprcdrifdef_quad2}  that, if we had totally ignored the stable component, i.e. if we had set $(I-\pi)\pjeps_t\Xeps=0$ at the very beginning of the analysis, 
we would miss the term $b_H^{q,(2)}$.
\end{rmk}

\begin{rmk}\label{rmk:explictcentmanterms_form}
The coefficients $b_H^{q,(k)}$ can be written more explicitly as 
\begin{align}\notag
b_H^{q,(1)}(\hlev)&=\frac{1}{2\pi/\om_c}\int_0^{2\pi/\om_c}dt\int_0^{(2\pi/\om_c)-t}ds\,\left( 2(\Psiz G_q(\etavg_t))^*\left[\begin{array}{cc}0 & e^{i\om_c s}\\ e^{-i\om_c s} & 0\end{array}\right]\Psiz G_{q}(\etavg_{t+s})\right) \\
& \qquad + \frac{\sqrt{2\hlev}}{2\pi/\om_c}\int_0^{2\pi/\om_c}dt\int_0^{(2\pi/\om_c)-t}ds\,\left((\Phi e^{sB}\Psiz G_q(\etavg_t)).\nabla\right)(\mathcal{E}_{t+s}G_q(\etavg_{t+s})), \label{eqn:explictcentmanterms_form_1}
\end{align}
\begin{align}
b_H^{q,(2)}(\hlev)&=\frac{\sqrt{2\hlev}}{2\pi/\om_c}\int_0^{2\pi/\om_c}dt\int_0^{\infty}ds\, \sum_{j=1}^n(G_q(\etavg_t))_j\left((\Th(s)(I-\pi)\Ind \ubar{e_j}).\nabla\right)(\mathcal{E}_{t+s}G_q(\etavg_{t+s})), \label{eqn:explictcentmanterms_form_2}
\end{align}
where $\etavg_t$ is defined in \eqref{eq:unpertsolused4avg}, and 
\begin{align}\label{eq:Ebyroot2h}
\mathcal{E}_t:=e^{-i\om_ct}\Psiz_1+e^{i\om_ct}\Psiz_2,
\end{align}
 and $ \ubar{e_j}$ denotes unit vector in the $j^{th}$ direction of $\R^n$.
To check how these explicit forms follow from \eqref{eq:taudef}--\eqref{eq:avgbhqk} refer to appendix \ref{subsec:PREpap:helpfulcomm}. If $G_q$ is a polynomial, the terms in \eqref{eqn:explictcentmanterms_form_1} can be put in Mathematica to get explicit functional dependence on $\hlev$; otherwise numerical integration can be done at specific $\hlev$ values. For the term in \eqref{eqn:explictcentmanterms_form_2} the integral $\int_0^{2\pi/\om_c}$ can be evaluated first using mathematica and then $\int_0^\infty$ can be done using numerical integration. All that we would need is the solutions of the unperturbed system with $n$ different initial conditions $(I-\pi)\Ind \ubar{e_j}$ for $j=1,\ldots,n$. Since the initial condition $(I-\pi)\Ind \ubar{e_j}$ belong to the stable space $Q$, the solution $\Th(s)(I-\pi)\Ind \ubar{e_j}$ decays exponentially fast to zero and hence then integral $\int_0^\infty$ need not be evaluated until infinity---a reasonable large upper limit would be enough to get a good enough approximation. An example is done next section to illustrate the above computations. Note that, when applied in a deterministic DDE setting, the above formulas provide an alternate way to compute the effect of center-manifold terms on the amplitude of critical mode.
\end{rmk}


\subsection{Example}\label{subsec:quadexampleCalcIllus} 
Consider the equation \eqref{eq:examplesys} with added quadratic nonlinearity $G_q(\eta)=(\eta(-1))^2$:
\begin{align}\label{eq:examplesys_quadadded}
dx(t)&=-\frac{\pi}{2}x(t-1)dt  + \eps^2 x^3(t-1)dt + \eps \sigma dW   + \eps x^2(t-1)dt
\end{align}
We apply theorem \ref{thm:PREpap:convgHforWquad}. Note that $b_H$ and $\sigma_H$ are already evaluated (see equations \eqref{eq:avgbHforscalarex} and \eqref{eq:avgsigHforscalarex}). We continue using the $\Phi$ and $\Psi$ from section \ref{sec:subsec:a_scalar_eq}.

Now we evaluate $b_H^{q,(1)}$ and $b_H^{q,(2)}$ using \eqref{eq:avgbhqk}. In section \ref{subsec:numverif} we show by numerical simulations how the averaged dynamics would be useful to gain information about \eqref{eq:examplesys_quadadded}.

Note that $(\xi.\nabla)G_q(\eta)=2\eta(-1)\xi(-1)$. We also write it as $2\eta\big|_{-1}\xi\big|_{-1}$ to avoid writing too many braces. Using the formula \eqref{eqn:explictcentmanterms_form_1}, we have $b_H^{q,(1)}(\hlev)=\frac{1}{2\pi/\om_c}\int_0^{2\pi/\om_c}\left(\int_0^{(2\pi/\om_c)-t}\mathscr{G}(t,s)\,ds\right)dt$ where
$$\mathscr{G}(t,s)=2\Psiz_1\Psiz_2(e^{i\om_c s}+e^{-i\om_c s})(\etavg_t\big|_{-1})^2(\etavg_{t+s}\big|_{-1})^2\,+\,\sqrt{2\hlev}\mathcal{E}_{t+s}2(\etavg_{t+s}\big|_{-1})(\Phi\big|_{-1}e^{sB}\Psiz)(\etavg_{t}\big|_{-1})^2,$$
where $\etavg_t$ is defined in \eqref{eq:unpertsolused4avg}.
Using Mathematica we get $b_H^{q,(1)}(\hlev)=-64\hlev^2/(4+\pi^2)^2 \approx -0.3327\hlev^2$.

To evaluate $b_H^{q,(2)}(\hlev)$ using \eqref{eqn:explictcentmanterms_form_2}, we first evaluate the $\int_0^{2\pi/\om_c}$ integral. We have
\begin{align*}
b_H^{q,(2)}(\hlev)&=\int_0^{\infty}\left(\frac{1}{2\pi/\om_c}\int_0^{2\pi/\om_c}\sqrt{2\hlev}\mathcal{E}_{t+s}2(\etavg_{t+s}\big|_{-1})(\Th(s)(I-\pi)\Ind \big|_{-1})(\etavg_{t}\big|_{-1})^2\,dt\right)\,ds \\
&=-\frac{4\hlev^2}{4+\pi^2}\int_0^\infty\big(2\pi+\pi\cos(\pi s)+2\sin(\pi s)\big)\, (\Th(s)(I-\pi)\Ind \big|_{-1})\,ds.
\end{align*} 
The $\int_0^\infty$ integral can be evaluated numerically by simulating the unperturbed system with the initial condition $(I-\pi)\Ind$, i.e. $\Ind-\Phi\Psiz$. We get $b_H^{q,(2)}(\hlev) \approx -0.7893\hlev^2$.


\subsection{Verification by numerical simulations}\label{subsec:numverif}

This section illustrates the results of theorems \ref{thm:PREpap:convgHforWnoquad} and \ref{thm:PREpap:convgHforWquad} using numerical simulations and also shows how the averaged $\hlev$ process can be used to gain information about the original $x$ dynamics (recall remark \ref{rmk:whyHprcUseful}). For details of the numerical scheme see appendix \ref{appsec:numericalscheme}.

Consider
\begin{align}\label{eq:examplesys_quadadded_forsimu}
dx(t)&=-\frac{\pi}{2}x(t-1)dt  + \eps^2 \gamma_c x^3(t-1)dt + \eps \sigma dW  + \eps \gamma_q x^2(t-1)dt.
\end{align}

Draw a random sample of size $N_{samp}$ with $\hlev$ values $\{\hlev^0_i\}_{i=1}^{Nsamp}$. Simulate them according to 
\begin{align}\label{eq:examplesysavgd}
d\hlev(t)=(b_H+b_H^{q,(1)}+b_H^{q,(2)})(\hlev(t))dt + \sigma_H(\hlev(t)) dW,
\end{align}
for $0\leq t \leq T_{end}$, where $b_H$ and $\sigma_H$ are obtained from \eqref{eq:avgbHforscalarex}, \eqref{eq:avgsigHforscalarex}, and $b_H^{q,(i)}$ are obtained in section \ref{subsec:quadexampleCalcIllus}:
\begin{align}\label{eq:auxyaux_driftdiffevals}
(b_H+b_H^{q,(1)}+b_H^{q,(2)})(\hlev)&=2\Psiz_1\Psiz_2\sigma^2-\gamma_c\frac32(i(\Psiz_1-\Psiz_2))  \hlev^2  -\gamma_q^2(0.3327+0.7893)\hlev^2, \\ \notag
\sigma_H^{2}(\hlev)&=4\Psiz_1\Psiz_2\sigma^2 \hlev. 
\end{align}

Fix $\eps$. Simulate \eqref{eq:examplesys_quadadded_forsimu} for $0 \leq t \leq T_{end}/\eps^2$ using initial history $\{\sqrt{2\hlev^0_i} \cos(\om_c\bullet)\}_{i=1}^{Nsamp}$. 

Fix a number $H^*$ and let $\tau^\eps$ be the first time $|x(t)|$ exceeds $\sqrt{2H^*}$ and $\tau^\hlev$ be the first time $\hlev(t)$ exceeds $H^*$, i.e.
\begin{align*}
\tau^\eps&:=\inf \{t\geq 0: |x(t)|\geq \sqrt{2H^*}\}, \\
\tau^\hlev&:=\inf \{t\geq 0: \hlev(t)\geq H^*\}.
\end{align*}

We can check whether the following pairs are close.
\begin{enumerate}
\item the distribution of $\ham(\pj_{T_{end}/\eps^2}x)$  from \eqref{eq:examplesys_quadadded_forsimu} (where $\ham$ is defined in \eqref{eq:hamfuncdef}) \emph{and} the distribution of $\hlev(T_{end})$ from \eqref{eq:examplesysavgd},
\item the distribution of $\eps^2 \tau^\eps$ \emph{and} the distribution of $\tau^\hlev$.
\end{enumerate}

\begin{figure}[h!]
    \centering
    \begin{minipage}{.5\textwidth}
        \centering
        \includegraphics[scale=0.4]{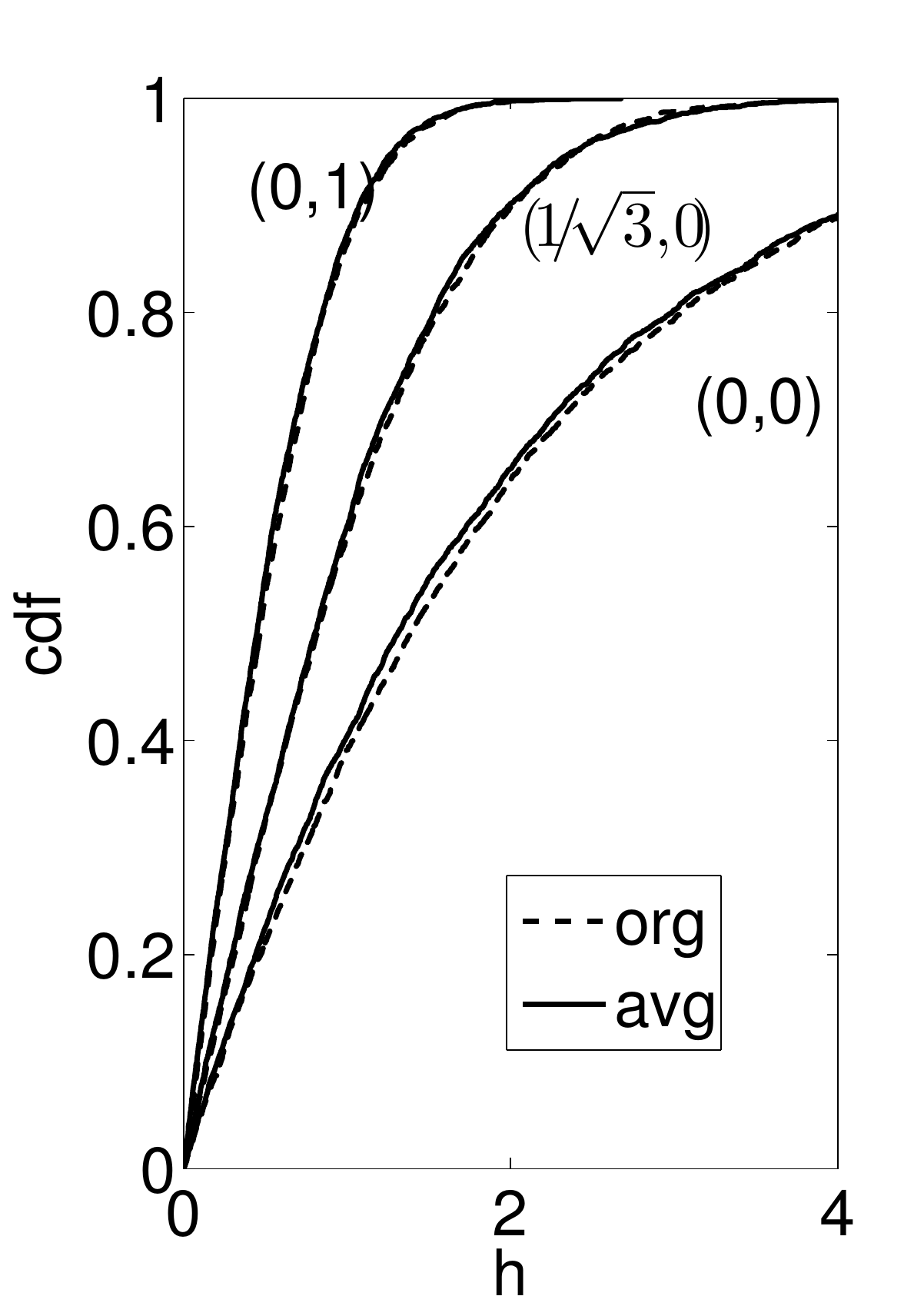}
        \caption{Cumulative distribution function (cdf) of $\ham(\pj_{2/\eps^2}x)$ ({\tt{org}}) and $\hlev(2)$ ({\tt{avg}}). The numbers in brackets are $(\gamma_q,\gamma_c)$ values.}
        \label{fig:tcdf}
    \end{minipage}%
    \begin{minipage}{0.5\textwidth}
        \centering
        \includegraphics[scale=0.4]{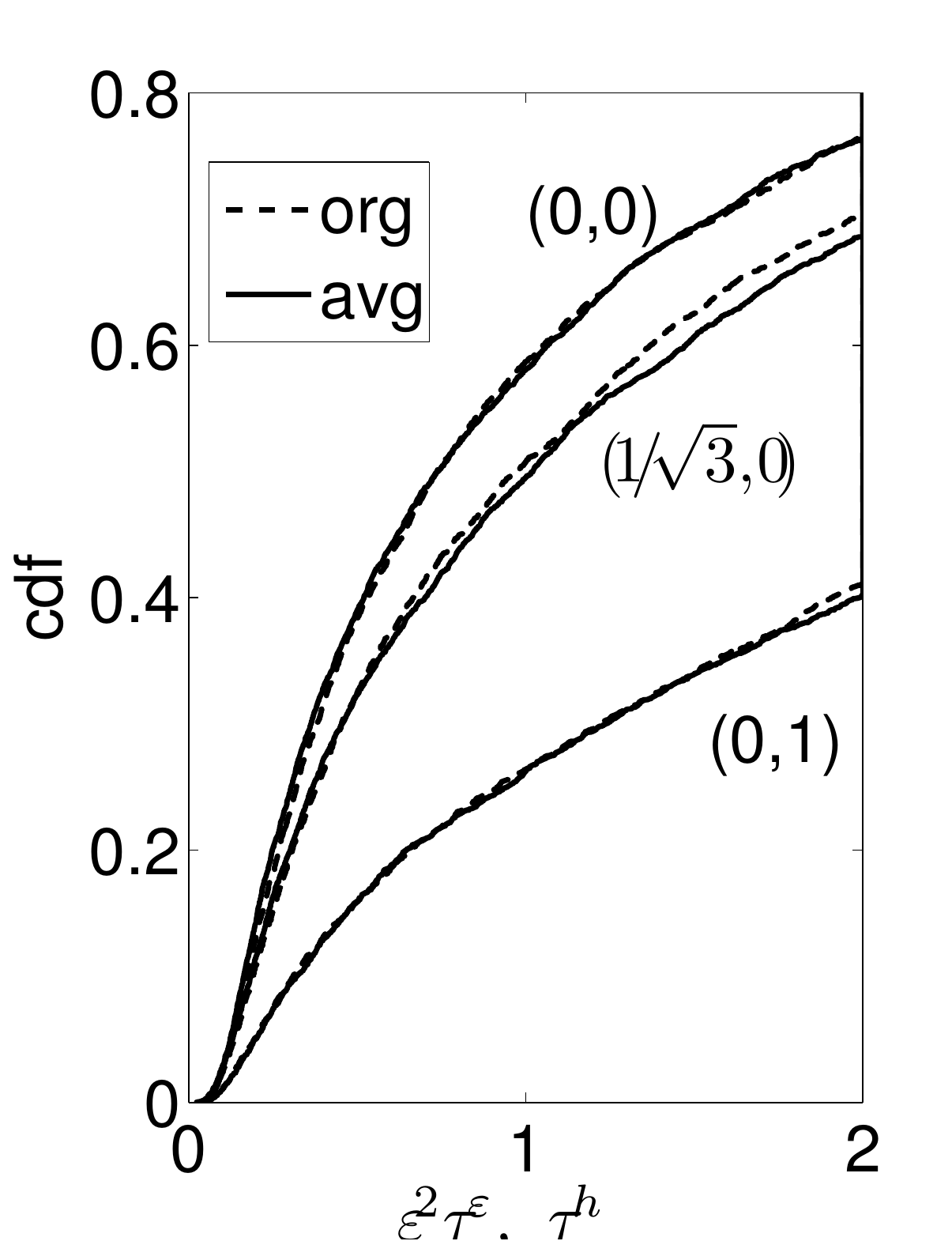}
         \caption{Cumulative distribution function (cdf) of $\eps^2 \tau^\eps$ ({\tt{org}}) and cdf of $\tau^\hlev$  ({\tt{avg}}). The numbers in brackets are $(\gamma_q,\gamma_c)$ values. The cdf value at $\eps^2\tau^\eps=2$ indicates the fraction of the sample whose modulus exceeded $\sqrt{2H^*}$ before the time $2/\eps^2$.}
        \label{fig:taucdf}
    \end{minipage}
\end{figure}

We took $\eps=0.025$, $H^*=1.5$, $T_{end}=2$, $N_{samp}=4000$, and $\sqrt{2\{\hlev^0_i\}_{i=1}^{Nsamp}}=1.2$.  Figures \ref{fig:tcdf} and \ref{fig:taucdf} answer the above questions. Three cases are considered with $\sigma=1$ fixed: $(\gamma_q=0,\gamma_c=0)$, $(\gamma_q=0,\gamma_c=1)$, $(\gamma_q=1/\sqrt{3},\gamma_c=0)$.

{\emph{ From the figures we can see that it is enough to study the averaged equations for $\ham(\pj_t x)$ to get a good approximation of the behaviour of $x$. The distribution of $\ham(\pj_t x)$ (note that $\sqrt{2\ham}$ gives the amplitude of oscillations) is well predicted by the distribution of the averaged system $\hlev$; and the distribution of time taken by $x$ to exceed a threshold $\sqrt{2H^*}$ is well predicted by the time taken by the averaged process $\hlev$ to exceed $H^*$.
Because the averaged equations do not contain any delay, they are easier to analyse and simulate numerically.}}


\section{Other kinds of noise}\label{sec:OTHnoises} 
Here we consider equations of the form
\begin{align}\label{eq:detDDE_pert_gennoise_nonon}
\begin{cases}dx(t)=L_0(\pj_t x)dt + \eps \sigma(\gnoise_t)F(\pj_t x)dt, \quad t\geq 0, \\
\pj_0x =\icond \in \C,
\end{cases}
\end{align}
where $F: \C \to \R^n$ is Lipschitz, with atmost linear growth and three bounded derivatives; and $\gnoise$ is a noise process whose state space is denoted by $\nstsp$, and $\sigma:\nstsp\to \R$. 

We make the following assumptions on the noise $\gnoise$.
\begin{assumpt}\label{assmp:assmp_on_noise}
The noise $\gnoise$ is a $\nstsp$-valued time-homogenous Markov process with transition probability function, $\nu$, given by
 $$\nu(t,\gnoise,B)=\mbbP\{\gnoise_t\in B \,| \,\gnoise_0=\gnoise\}$$
for $B$ a borel subset of $\nstsp$. There exist a unique invariant probability measure $\bar{\nu}$ and positive constants $c_1$ and $c_2$ such that for all $t\geq 0$,
$$\sup_{\gnoise \in \nstsp}\int_{\nstsp}|\nu(t,\gnoise,d\zeta)-\bar{\nu}(d\zeta)| \leq c_1e^{-c_2t},$$
i.e. the transition probability density converges to stationary density exponentially fast.  
The function $\sigma$ is bounded,  and such that $\int_{\nstsp}\sigma(\gnoise)\bar{\nu}(d\gnoise)=0$. 

Other requirements are: $\nstsp$ is locally compact separable metric space; the transition semigroup is Feller with $\sigma(\cdot)$ in the domain of the infinitesimal generator.
\end{assumpt}

For example, a finite-state continuous-time markov chain satisfies the above requirements.

The autocorrelation of the noise process $\gnoise$ is denoted by $R$:
\begin{align}\label{eq:PREpap:autocorr}
R(s)=\int_{\nstsp} \sigma(\gnoise)\left(\int_{\nstsp}\sigma(\zeta)\,\nu(s,\gnoise,d\zeta)\right)\bar{\nu}(d\gnoise).
\end{align}

For the perturbed system \eqref{eq:detDDE_pert_gennoise_nonon}, $\ham(\pj_t x)$ varies slowly compared to $x$. Changes in $\ham(\pj_t x)$ are significant only on times of order $1/\eps^2$. Hence, we rescale time and write $X^\eps(t)=x(t/\eps^2)$ where $x$ is governed by \eqref{eq:detDDE_pert_gennoise_nonon}. Also, we write $\gnoiseeps_t=\gnoise(t/\eps^2)$.

Using the segment extractor $\pjeps_t$ defined in \eqref{eq:newsegextract}, $X^\eps$ satisfies
\begin{align}\label{eq:detDDE_pert_rescale_gennoise}
\begin{cases}dX^\eps(t)=\frac{1}{\eps^2}L_0(\pjeps_t \Xeps)dt+ \frac{1}{\eps}\sigma(\gnoiseeps_t) F(\pjeps_t \Xeps)dt, \quad t\geq 0, \\
\pjeps_0\Xeps =\icond \in \C.
\end{cases}
\end{align}

Write $\hprc^{\eps}(t):=\ham(\pjeps_t \Xeps)$. Then $\hprc^{\eps}(t)$ satisfies
\begin{align}\label{eq:evolofhprceps_gennoise}
d\hprc^{\eps}(t)=\frac{1}{\eps}\sigma(\gnoiseeps_t)\hprcdrift(\pjeps_t \Xeps)dt, \qquad \qquad \hprc^{\eps}(0)=\ham(\icond)
\end{align}
where 
\begin{align}
\hprcdrift(\eta)&=E(\eta) F(\eta), \label{eq:hprcdiffdef_gennoise} 
\end{align}
where $E$ is defined in \eqref{eq:hprcdrifE}. 

Using the technique of martingale problem, we can prove\footnote{Proof of theorem \eqref{thm:gennoiseweakconvg} and a result to the effect that the stable component of the solution is small would be published in a different article.} the following result (a sketch of proof is given in appendix \ref{appsec:sketchproofGenNoise}):
\begin{thm}\label{thm:gennoiseweakconvg}
Under the conditions on $F$ and noise $\gnoise$ listed before; the probability distribution of $\hprc^{\eps}$ converges, as $\eps \to 0$,  to the distribution of the process $\avgHproc$ which is the solution of the SDE
$$d\avgHproc(t)=b_H(\avgHproc(t))dt+\sigma_H(\avgHproc(t))dW(t), \qquad \quad \avgHproc(0)=\ham(\icond),$$
with coefficients $b_H$ and $\sigma_H$ given by
\begin{align*}
\sigma_{H}^2(\hlev)&=\frac{1}{2\pi/\om_c}\int_0^{2\pi/\om_c}2\, \hprcdrift(\etavg_t) \,\left(\int_0^\infty R(s)\,\hprcdrift(\etavg_{t+s})\,ds\right)dt,\\
b_H(\hlev)&=\frac{1}{2\pi/\om_c}\int_0^{2\pi/\om_c}\left(\int_0^\infty R(s)\, \left(\Th(s)\Ind F(\etavg_t).\nabla\right)\,\hprcdrift(\etavg_{t+s}) \,ds\right)dt,
\end{align*}
where $\etavg_t$ is defined in \eqref{eq:unpertsolused4avg}.
\end{thm}
We urge the reader to study appendix \ref{appsec:sketchproofGenNoise} to gain intuition about the process of obtaining the coefficients $b_H$ and $\sigma_H$. Akin to the formulas \eqref{eqn:explictcentmanterms_form_1}--\eqref{eqn:explictcentmanterms_form_2}, the coefficient $b_H$ can be written more explicitly as
\begin{align*}\notag
b_H(\hlev)&=\frac{1}{2\pi/\om_c}\int_0^{2\pi/\om_c}dt\int_0^{\infty}ds\,\left( 2R(s)\,(\Psiz F(\etavg_t))^*\left[\begin{array}{cc}0 & e^{i\om_c s}\\ e^{-i\om_c s} & 0\end{array}\right]\Psiz F(\etavg_{t+s})\right) \\
& \qquad + \frac{\sqrt{2\hlev}}{2\pi/\om_c}\int_0^{2\pi/\om_c}dt\int_0^{\infty}ds\,R(s)\sum_{j=1}^n(F(\etavg_t))_j\left((\Th(s)\Ind \ubar{e_j}).\nabla\right)(\mathcal{E}_{t+s}F(\etavg_{t+s})), 
\end{align*}
where $\etavg_t$ is defined in \eqref{eq:unpertsolused4avg}, $\mathcal{E}$ is defined in \eqref{eq:Ebyroot2h}, and $\ubar{e_j}$ is the unit vector in the $j^{th}$ direction of $\R^n$. Similarly,
\begin{align*}
\sigma_{H}^2(\hlev)&=\frac{4\hlev}{2\pi/\om_c}\int_0^{2\pi/\om_c}dt\int_0^{\infty}ds\, (\mathcal{E}_{t}F(\etavg_t)) R(s) (\mathcal{E}_{t+s}F(\etavg_{t+s}))).
\end{align*}
It would be easier to do the $\int_0^{2\pi/\om_c}$ integral before the $\int_0^\infty$ integral.

Analogous results for systems without delay are found in section 4 of \cite{PapKoh75}. Even systems with delay can be put in the framework of \cite{PapKoh75}. Equations of the form \eqref{eq:detDDE_pert_gennoise_nonon} with $F(0)=0$ and $\int_{\nstsp}\sigma(\gnoise)\bar{\nu}(d\gnoise)\neq 0$ (i.e noise is not mean zero) are studied in \cite{Tsark1}.

\begin{rmk}
In the equation \eqref{eq:detDDE_pert_gennoise_nonon}, we could have included the deterministic perturbations $G$ and $G_q$ as done in equation \eqref{eq:quadnon_main_in_short_form}; but the averaged drift terms arising from these would be same as in the previous sections.
\end{rmk}


\subsection{Linear perturbations}\label{subsec:realnoiselinpert_v} 
When $F(\eta)=L_1\eta$ where $L_1:\C \to \R^n$ is a linear operator, the expressions for $b_H$ and $\sigma_H$ can be more explicitly evaluated using the autocorrelation function as follows. Let $\Upsilon$ be the $2\times 2$ matrix $\Upsilon_{ij}=\Psiz_i L_1\Phi_j$. Let 
\begin{align*}
R_0&=\int_0^\infty R(s)ds, \\
R_{2c}&=\int_0^\infty R(s)\cos(2\om_c s)ds,\\
\hat{R}_1&=\int_0^\infty R(s)e^{-i\om_c s}\Psiz_1L_1(\Th(s)(I-\pi)\Ind L_1\Phi_1)\,ds, \\
\hat{R}_2&=\int_0^\infty R(s)e^{i\om_c s}\Psiz_2L_1(\Th(s)(I-\pi)\Ind L_1\Phi_2)\,ds.
\end{align*}
Then, 
$$b_H(\hlev)=C_b\hlev, \qquad \sigma_H^2(\hlev)=C_{\sigma}\hlev^2$$
where
\begin{align*}
C_b&=\bigg((\Upsilon_{11}+\Upsilon_{22})^2R_0+4\Upsilon_{12}\Upsilon_{21}R_{2c}+\hat{R}_1+\hat{R}_2\bigg), \\
C_{\sigma}&=2\bigg((\Upsilon_{11}+\Upsilon_{22})^2R_0+2\Upsilon_{12}\Upsilon_{21}R_{2c}\bigg).
\end{align*}

\begin{rmk}\label{rmk:Foferror2}
Note that if we had totally ignored the stable modes, i.e. if we set $(I-\pi)\pjeps_t\Xeps=0$ at the very beginning of the analysis, 
we would not have the terms $\hat{R}_1$ and $\hat{R}_2$.
\end{rmk}

The Lyapunov exponent for the averaged equation
\begin{align}\label{mult:eq:examplesysavgd_gennoiselinear}
d\hlev(t)\,=b_H(\hlev) dt\,+\,\sigma_H(\hlev)\,dW,
\end{align}
can be calculated to be 
\begin{align}\label{eq:linpert_gennoise_lamavg}
\lambda_{avg}=C_b-\frac12C_{\sigma}=2\Upsilon_{12}\Upsilon_{21}R_{2c}+\hat{R}_1+\hat{R}_2.
\end{align}

Using singular perturbation methods and Furstenberg-Khasminskii formula, the following theorem for \emph{scalar} processes is proved in \cite{NavWihs} and \cite{WihsNavam}.
\begin{thm}\label{thm:gennoiselyap}
Consider \eqref{eq:detDDE_pert_gennoise_nonon} with $F(\eta)=L_1(\eta)$ where $L_1:\C \to \R$ is linear. Let the top Lyapunov exponent of the process $x$ be defined by
\begin{align}\label{eq:topLyapexpDef}
\lambda^\eps:=\limsup_{t\to \infty}\frac1t\ln \sup_{s\in [t-r,t]}|x(s)|.
\end{align}
 Then $\lambda^\eps  = \eps^2\frac12\lambda_{avg}+O(\eps^3)$.
\end{thm}
The same can be said about vector valued processes.


\subsection{Verification by numerical simulation}\label{subsec:verifyNumSimMarkovChain} 
Consider the system
\begin{align}\label{eq:detDDE_pert_gennoise_numsim}
dx(t)=-\frac{\pi}{2}x(t-1)dt + \eps \sigma(\gnoise_t)x(t-1)dt.
\end{align}

Let $\gnoise$ be a two-state symmetric markov chain with switching rate $g/2$, i.e. 
\begin{align}\label{eq:rateofswitch2state}
\lim_{t\downarrow 0}\frac1t P_{1\to 2}(t)=g/2=\lim_{t\downarrow 0}\frac1t P_{2\to 1}(t)
\end{align}
 where $P_{i\to j}(t)$ is the probability of transition from state $i$ to state $j$ in time $t$. Let $\sigma(\gnoise=1)=-\sigma(\gnoise=2)=\sigma_0$. We then have the autocorrelation as $R(s)=\sigma_0^2e^{-gs}$. 

We consider two cases $g=2$ or $g=6$ with $\sigma_0=1$. The averaged equations are
\begin{align*}
g=2: \qquad &d\hlev(t)\,=0.3734\,\hlev\, dt\,+\,\sqrt{0.9873}\,\hlev\,dW,\\
g=6: \qquad &d\hlev(t)\,=0.1715\,\hlev\, dt\,+\,\sqrt{0.4245}\,\hlev\,dW.
\end{align*}

 Using same notation as in section \ref{subsec:numverif}, we fix $\eps=0.025$, $T_{end}=1$, $H^*=1$, $N_{samp}=4000$ and $\sqrt{2\{\hlev^0_i\}_{i=1}^{Nsamp}}=1$. The equation \eqref{eq:detDDE_pert_gennoise_numsim} is simulated for time $T_{end}/\eps^2$ with initial history $\{\sqrt{2\hlev^0_i} \cos(\om_c\bullet)\}_{i=1}^{Nsamp}$.  We obtain the following figures \ref{fig:tcdf_twostate} and \ref{fig:taucdf_twostate}
\begin{figure}[h!]
    \centering
    \begin{minipage}{.5\textwidth}
        \centering
        \includegraphics[scale=0.45]{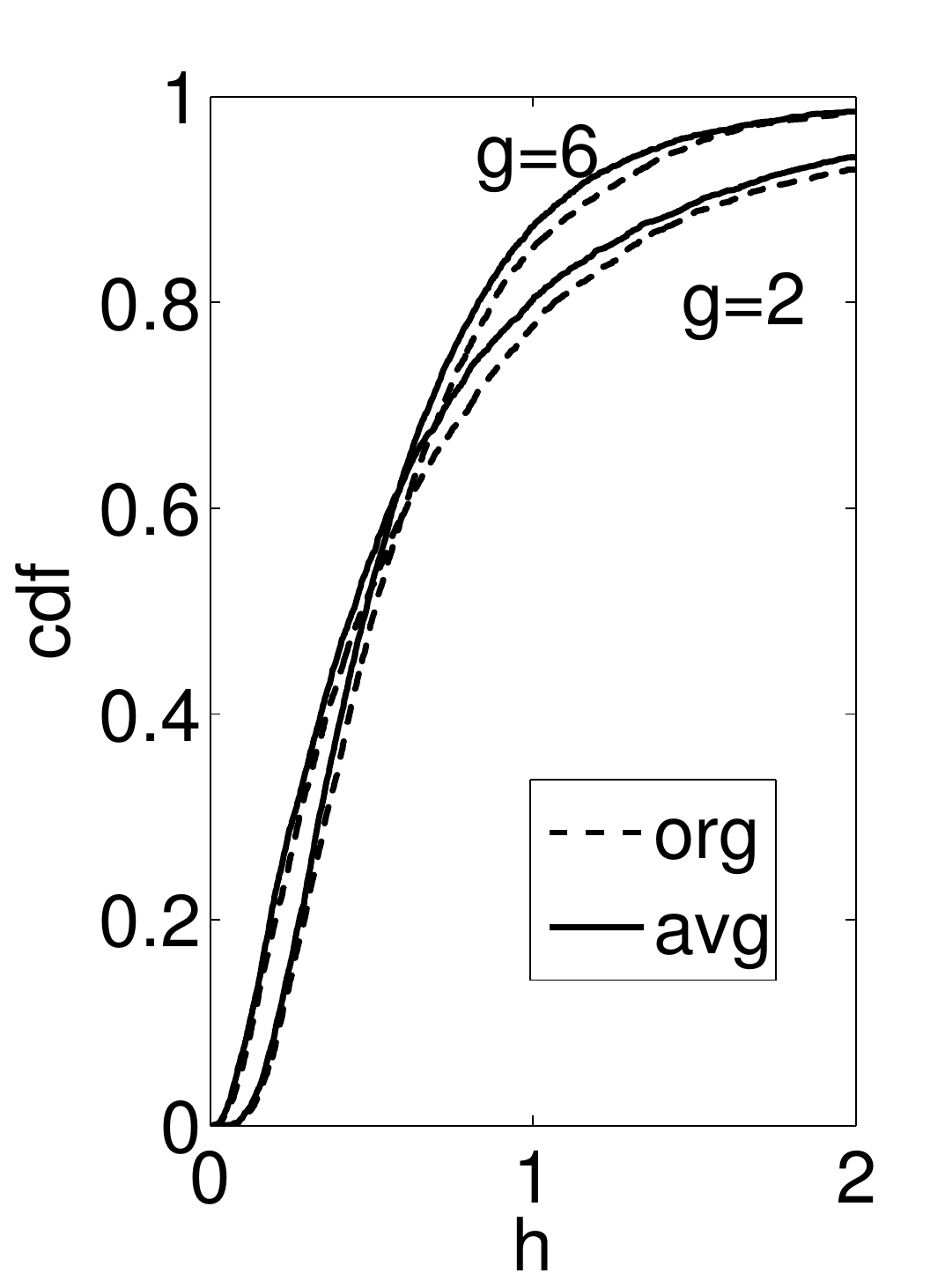}
        \caption{Cumulative distribution function (cdf) of $\ham(\pj_{1/\eps^2}x)$ ({\tt{org}}) and $\hlev(1)$ ({\tt{avg}}).}
        \label{fig:tcdf_twostate}
    \end{minipage}%
    \begin{minipage}{0.5\textwidth}
        \centering
        \includegraphics[scale=0.4]{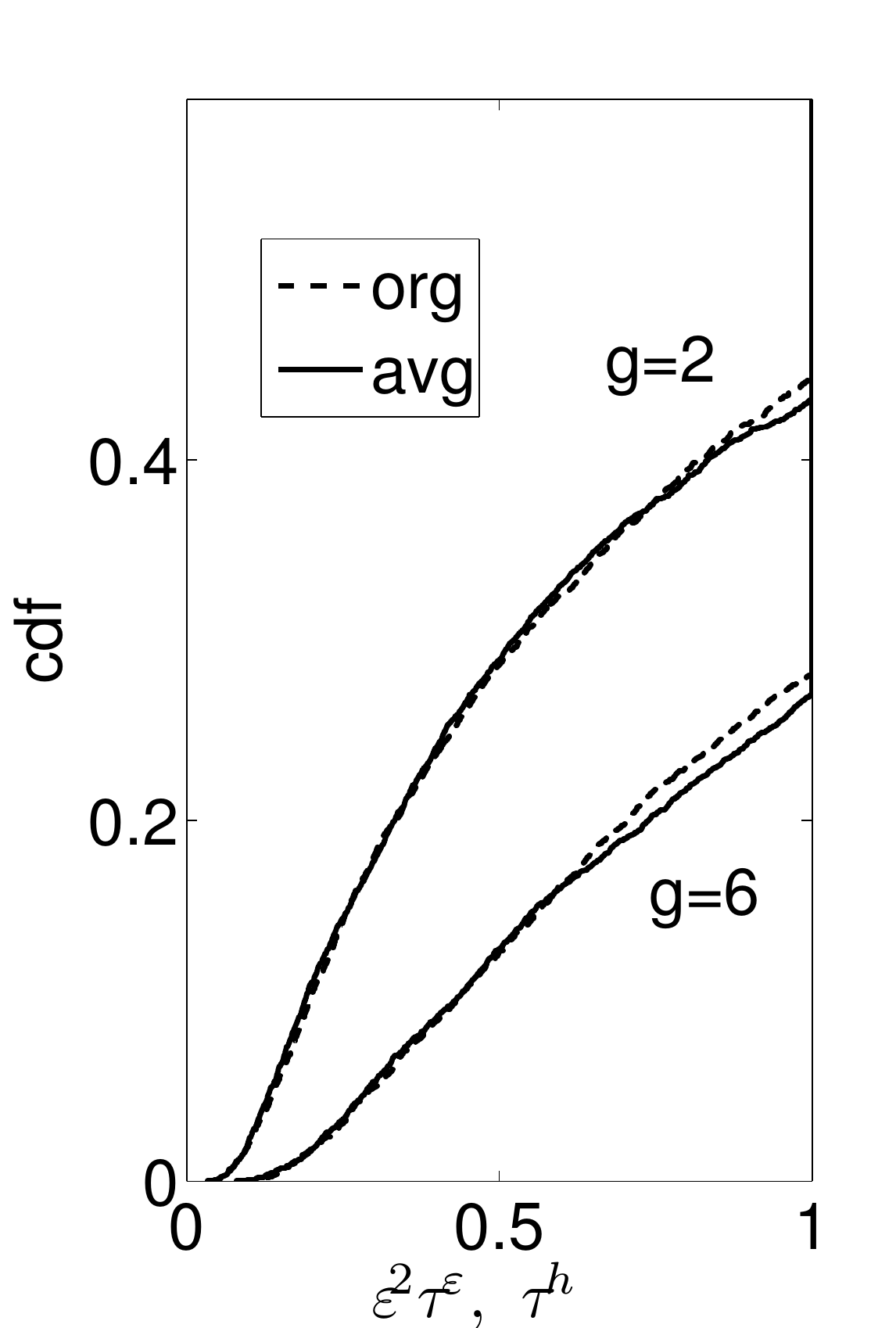}
        \caption{Cumulative distribution function (cdf) of $\eps^2 \tau^\eps$ ({\tt{org}}) and cdf of $\tau^\hlev$  ({\tt{avg}}). The cdf value at $\eps^2\tau^\eps=1$ indicates the fraction of particles whose modulus exceeded $\sqrt{2H^*}$ before the time $1/\eps^2$.}
        \label{fig:taucdf_twostate}
    \end{minipage}
\end{figure}
which show that the averaged system gives a good approximation of the original system. For details of the numerical scheme see appendix \ref{appsec:numericalscheme}.


\section{Discussion}\label{sec:Discussion} 

Delay equations with noise perturbations as considered in section \ref{sec:OTHnoises} display interesting similarities with non-delay systems. For example, \cite{NSNVedFIOsc} considers coupled oscillators with one of the oscillators stable, in the following form. Let $J$ be the symplectic matrix $\left(\begin{array}{cc}0&1\\-1&0\end{array}\right)$, $I$ be the $2\times 2$ identity matrix and $O$ be the $2\times 2$ zero matrix. Let $x\in \R^4$ be governed by
\begin{align}\label{eq:NavVed_couposcstab}
\dot{x}(t)=\left(\begin{array}{cc}\om_1 J  & O\\ O & -\delta I+\om_2J \end{array}\right)x(t) + \eps \sigma(\gnoise(t))\left(\begin{array}{cc}K & M\\ N & L \end{array}\right)x(t)
\end{align}
where $K, L, M, N$ are $2\times 2$ matrices. The oscillator with frequency $\om_1$ is coupled to the stable oscillator of frequency $\om_2$. \cite{NSNVedFIOsc} shows that the Lyapunov exponent of the above system can be written in terms of quantities analogous to $R_0,\,R_{2c},\,\hat{R}_i$ defined in section \ref{subsec:realnoiselinpert_v}. Further they show that both stabilization and destabilization are possible depending on the matrix coefficients $K, M$ and $N$.

\begin{figure}
\centering
\includegraphics[scale=0.65]{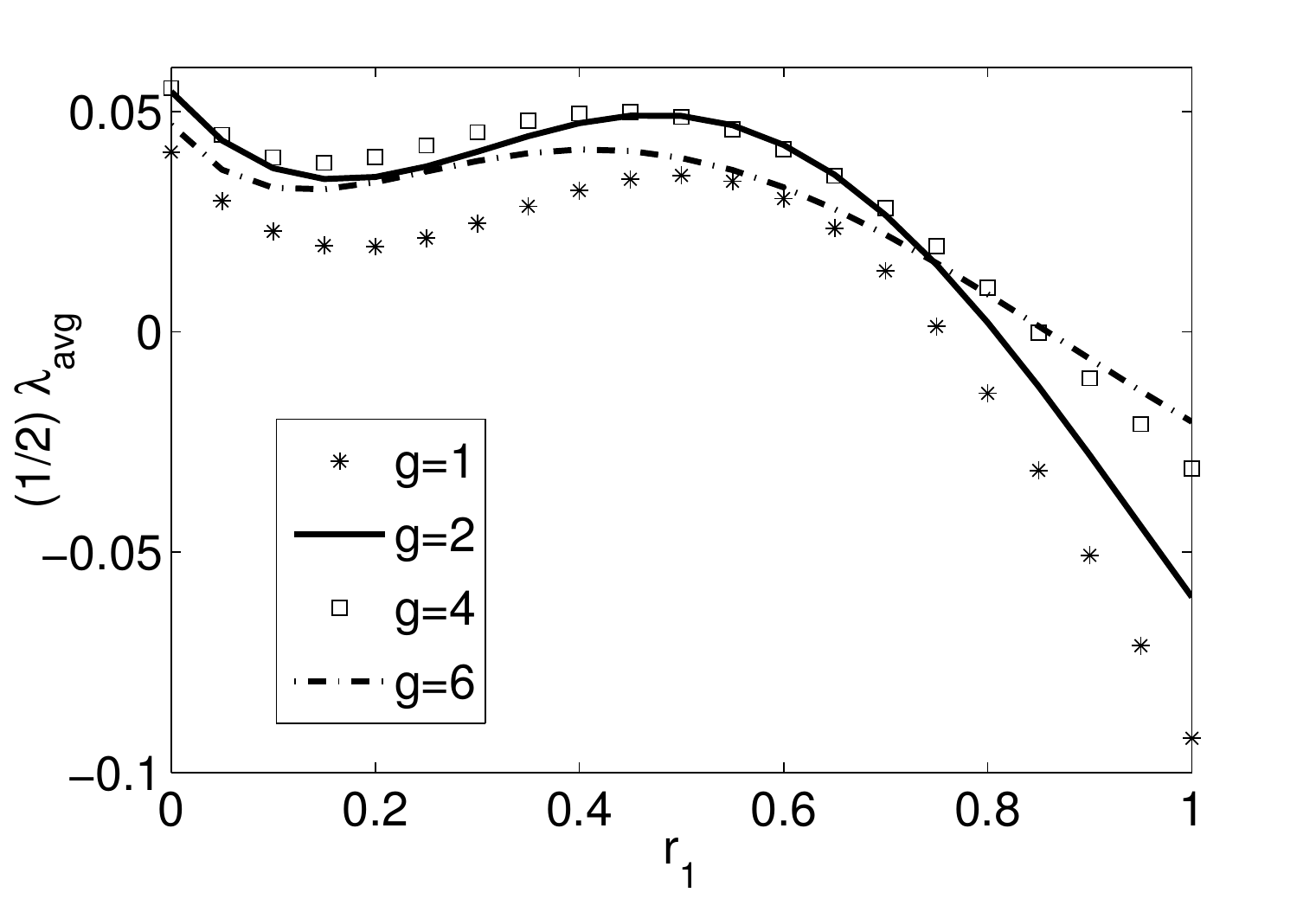}
  \caption{$\frac12\lambda_{avg}$ as a function of the delay in the perturbation ($r_1$) and the rate of switching of the noise ($g$) for equation \eqref{eq:LyapVsDel_illus}. The top Lyapunov exponent $\lambda^\eps$ is close to $\eps^2\frac12\lambda_{avg}$ by theorem \ref{thm:gennoiselyap}. Note that both $\lambda_{avg}<0$ (stabilization) and $\lambda_{avg}>0$ (destabilization) are possible.}
  \label{fig:LyapVsDel}
\end{figure}

The delay system that we considered under the assumption \ref{ass:assumptondetsys} can be thought of as a coupled oscillator system with one critical mode and infinitely many stable modes (the characteristic equation has a pair of roots $\pm i\om_c$, and all other roots have negative real part).  The lyapunov exponent obtained in \eqref{eq:linpert_gennoise_lamavg} suggests that both stabilization and destabilization are possible. To illustrate this, consider 
\begin{align}\label{eq:LyapVsDel_illus}
dx(t)=-\frac{\pi}{2}x(t-1)dt + \eps \sigma(\gnoise_t)x(t-r_1)dt
\end{align}
with $\gnoise$ a two-state symmetric markov chain with states $\sigma(\gnoise)\in \{+1,-1\}$ and rate of switching $g/2$ (defined in \eqref{eq:rateofswitch2state}). Theorem \ref{thm:gennoiselyap} says that the Lyapunov exponent $\lambda^\eps$ (defined in \eqref{eq:topLyapexpDef}) is close to $\eps^2\frac12\lambda_{avg}$ where $\lambda_{avg}$ is evaluated in \eqref{eq:linpert_gennoise_lamavg}. Figure \ref{fig:LyapVsDel} shows how  $\frac12\lambda_{avg}$ varies with the delay in the perturbation ($r_1$) and rate of switching ($g$) of the two-state markov chain. Note that both $\lambda_{avg}<0$ (stabilization) and $\lambda_{avg}>0$ (destabilization) are possible.

Even the white noise allows for both possibilites. As mentioned in section \ref{subsec:linearpertWhitenoise}, the lyapunov exponent $\lambda_{avg}$ corresponding to \eqref{mult:eq:examplesys} equals $-\,Re[(\Psiz_1 L_1\Phi_1)^2]$. Applying to $dx(t)=-\frac{\pi}{2}x(t-1)dt + \eps x(t-r_1)dW$ we find that $\lambda_{avg}<0$ for $r_1<0.8609$ and $\lambda_{avg}>0$ for $0.8609<r_1\leq 1$.

The above examples raise the question whether stabilization or destabilization is possible when the noise is additive, i.e. the coefficient $F$ is a constant independent of the state $x$. To answer this question consider
\begin{align}\label{eq:additivenoiseHutt_eqn}
d\tilde{x}(t)&=\left(-\frac{\pi}{2} - \eps^2\gamma_o\right)\tilde{x}(t-1)dt  + \gamma_q \tilde{x}^2(t-1)dt +\gamma_c \tilde{x}^3(t-1)dt + \eps^2 \sigma dW.  
\end{align}
Scaling according to $\tilde{x}(t/\eps^2)=\eps X^\eps(t)$ we find that $X^\eps$ has same distribution as equation \eqref{eq:detDDE_pert_rescale_quadadded} with $L_0\eta = -\frac{\pi}{2}\eta(-1)$, $G_q(\eta)=\gamma_q\eta^2(-1)$, $G(\eta)=\gamma_c\eta^3(-1)-\gamma_o\eta(-1)$ and $F(\eta)=\sigma$.
The averaged equation corresponding to this is (obtained by evaluation of quantities in \eqref{eq:auxyaux_driftdiffevals} of section \ref{subsec:numverif} using $\Psiz_i$ from section \ref{sec:subsec:a_scalar_eq})
\begin{align*}
d\hlev(t)=\mathfrak{B}(\hlev(t))dt + \sigma  \sqrt{2\times 0.5768\,\hlev(t)} \,dW,
\end{align*}
where 
\begin{align}\label{eq:additivenoiseHutt_drift}
\mathfrak{B}(\hlev)=0.5768\sigma^2+0.9060\gamma_o\hlev-(1.3591\gamma_c+1.1220\gamma_q^2)\hlev^2.
\end{align} 
Let $(\gamma_q,\gamma_c)$ be such that $\widehat{\gamma}\overset{\text{def}}=1.3591\gamma_c+1.1220\gamma_q^2>0$. Assume that the noise is absent, i.e. $\sigma=0$. If $\gamma_o<0$ then $\hlev=0$ is the only fixed point\footnote{Fixed points are obtained by solving $\mathfrak{B}(h)=0$.}  and it is stable. If $\gamma_o>0$ then the zero fixed point looses stability and another stable fixed point $\hlev=0.9060\gamma_o/\widehat{\gamma}$ exists. \emph{In the presence of noise} ($\sigma\neq 0$), irrespective of the sign of $\gamma_0$, there are no fixed points because the diffusion is non-zero everywhere except at zero, and at zero $\mathfrak{B}(0)\neq 0$. Thus the additive noise destroys the fixed points. The amplitude ($A$) of oscillations is approximately $\sqrt{2\hprc}$ (recall remark \ref{rmk:whyHprcUseful}). The averaged equation corresponding to the amplitude $A=\sqrt{2\hprc}$ is (applying Ito formula), $dA=\frac{1}{A}(\mathfrak{B}(A^2/2)-\frac{c}{2}\sigma^2)dt+\sigma\sqrt{c}\,dW$ where $c=0.5768$. 

\cite{HuttEPLosc} considers stability of scalar delay systems with additive white noise. \cite{HuttEPLosc} writes equations for the individual projections $\la \Psi_i,\pjeps_t\Xeps \ra$, and using formal higher order corrections to the center-manifold, arrive at a differential equation for the mean of the amplitude of oscillations. However \cite{HuttEPLosc} commits the error of taking the mean of individual projections $\la \Psi_i,\pjeps_t\Xeps \ra$ to arrive at the mean of amplitude. The correct way to do is to take the mean of $\sqrt{2\hprc}$, i.e., $2\sqrt{\la \Psi_1,\pjeps_t\Xeps \ra \la \Psi_2,\pjeps_t\Xeps \ra}$. Though higher order corrections are provided by \cite{HuttEPLosc}, this error would have the effect of dropping of the term $0.5768\sigma^2$ in \eqref{eq:additivenoiseHutt_drift}. To understand the nature of the error more clearly, one can ignore all nonlinearities and consider $dx=-\frac{\pi}{2}x(t-1)dt+ \eps\sigma dW$. For this equation, the analysis in \cite{HuttEPLosc} predicts that the mean of the amplitude of oscillations does not change at all. However, the averaging results of this article predicts that $h$ evolves according to $dh=c\sigma^2dt+\sigma \sqrt{2ch}dW$ where $c=0.5768$ and hence the amplitude $A=\sqrt{2h}$ evolves according to (applying Ito formula) $dA=\frac12\sigma^2c\frac{1}{A}dt+\sigma \sqrt{c}dW$, from which we get that the mean of the amplitdue changes with time. \cite{HuttEPLosc} arrives at the conclusion that additive noise has the ability to postpone (stabilize) the bifurcation. However, in light of the above mistake\footnote{\cite{HuttEPLosc} also commits one more error of the following nature. In passing from equation 12 to equation 13 in \cite{HuttEPLosc} they assume that $\expt[g(X)]=g(\expt X)$ for a nonlinear function $g$ of a random variable $X$ (here $\expt$ denotes the expected value). This is wrong. For example, from the SDE, $dA_t=g(A_t)dt+\sigma dW_t$, one cannot claim that the mean varies as $\frac{d}{dt}\expt A_t = g(\expt A_t)$. In general, the moments of $A$ of different orders are coupled by the nonlinearity in the drift.} the conclusion must be re-evaluated.

The averaging results presented in this article allow us to simplify the analysis of delay systems at the verge of instability. The averaged dynamics does not involve any delay and hence is easier to analyse. Using numerical simulations we have amply demonstrated the usefulness of the theoretical results in approximating the probability distribution of the time-delay system with that of the averaged system. In section \ref{subsec:vdpolosc} we have shown how these results would be useful in computing an approximation to the shift of bifurcation thresholds in presence of noise.

We conclude this article with a section on a different instability scenario.

\section{A different kind of instability}

The instability in assumption \ref{ass:assumptondetsys} is not the only kind of instability possible. For example, one can have
\begin{assumpt}\label{ass:assumptondetsys_Hutt}
The characteristic equation \eqref{eq:chareq} has zero as a simple root, and all other roots have negative real parts.
\end{assumpt}
The analysis under assumption \ref{ass:assumptondetsys_Hutt} is similar to the analysis in previous sections.
Choose $\ubar{d}$ such that $\Delta(0)\ubar{d}=0_{n\times 1}$ and $\ubar{d_2}$ such that $\ubar{d_2}\Delta(0)=0_{1\times n}$. Define $\Phi$ by the constant $\Phi(\bullet)=\ubar{d}$ and $\Psi$ by $\Psi(\bullet)=c\ubar{d_2}$ where the constant $c$ is choosen so that $\la \Psi,\Phi \ra =1$ for the bilinear form in \eqref{eq:bilinform}. The space $\C$ can be split as $\C=P\oplus Q$ where $P$ is the space spanned by the constant function $\Phi$. The projection operator is $\pi:\C \to P$ given by $\pi(\eta)=\Phi\la \Psi,\eta\ra$. Define $\Psiz \overset{\text{def}}=\Psi(0)$.
Let $\Th$ and $\Ind$ be as defined in section \ref{subsec:strongdetpert_not}.
For the unperturbed system \eqref{eq:detDDE}, writing $\pj_tx = \pi \pj_tx+ (1-\pi)\pj_tx =\Phi z(t)+(I-\pi)\pj_t x$ we find that $\dot{z}=0$ and $||(I-\pi)\pj_tx||$ decays exponentially fast. So, defining $\ham(\eta)=\la \Psi,\eta\ra$ we find that $\hprc(t)=\ham(\pj_tx)$ is a constant for the unperturbed system (note that $\hprc$ is same as $z$). Now consider equations of the form \eqref{eq:quadnon_main_in_short_form}. Akin to condition \eqref{eq:assumption_on_Gq_zero} we need to impose that 
\begin{align}\label{eq:assumption_on_Gq_zero_Hutt}
\Psiz G_q(\Phi h)=0, \quad \forall h\in \R.
\end{align}
(\emph{If the above is not imposed, then the distribution of $\hprc$ on times of order $1/\eps$ converges to that of a deterministic process given by $\dot{\hprc}=\Psiz G_q(\Phi\hprc)$. Remark \ref{rmk:whencenteringnotimposed} deals with the case when \eqref{eq:assumption_on_Gq_zero_Hutt} is not satisfied.}) When \eqref{eq:assumption_on_Gq_zero_Hutt} is imposed, significant changes in $\hprc$ occurs only on times of order $1/\eps^2$. So writing $\Xeps(t)=x(t/\eps^2)$ we find that $\Xeps$ has the same probability distribution as the process satisfying \eqref{eq:detDDE_pert_rescale_quadadded}. Defining $\hprc^\eps(t):=\ham(\pjeps_t\Xeps)$ and using Ito formula we get that $\hprc$ satisfies \eqref{eq:evolofhprceps_quad} with $\hprcdrift^{q,(1)}(\eta)=\Psiz G_q(\pi\eta)=0$, $\hprcdrift^{q,(2)}(\eta)=\Psiz(G_q(\eta)-G_q(\pi\eta))$, $\hprcdrift(\eta)=\Psiz G(\eta)$ and $\sigma(\eta)=\Psiz F(\eta)$. It can be shown that result analogous  to theorem \ref{thm:PREpap:convgHforWquad} holds with the averaged drift and diffusion coefficients given by $b_H(\hlev)=\Psiz G(\Phi h)$, $\sigma_H^2(\hlev)=(\Psiz F(\Phi \hlev))^2$, $b_H^{q,(1)}=0$, and 
\begin{align}\label{eq:HuttBifScenario_whiteNoise_cent}
b_H^{q,(2)}(\hlev)=\int_0^\infty ((\Th(s)(I-\pi)\Ind G_q(\Phi \hlev)).\nabla)\Psiz G_q(\Phi \hlev)ds.
\end{align}
For \emph{scalar} systems the condition \eqref{eq:assumption_on_Gq_zero_Hutt} would necessarily mean that $G_q(\Phi \hlev)=0$ which would result in $\Ind G_q(\Phi \hlev)=0$ and hence $b_H^{q,(2)}=0$. This means that, when \eqref{eq:quadnon_main_in_short_form} is scalar valued, $G_q$ terms would have negligible effect on the dynamics on $P$ subspace for times of order $1/\eps^2$.

\begin{rmk}\label{rmk:whencenteringnotimposed}
When \eqref{eq:assumption_on_Gq_zero_Hutt} is not satisfied, the distribution of $\hprc$ on times of order $1/\eps$ converges to that of a deterministic process given by $\dot{\hprc}=\Psiz G_q(\Phi\hprc)$. Stochastic limit can be obtained if the strength of the noise is increased from $\eps$ to $\sqrt{\eps}$. Consider
\begin{align}
dx(t)=L_0(\pj_t x)dt\,&+\,\eps G_q(\pj_t x)dt\,+\,\eps^2 G(\pj_t x)dt   +\,\sqrt{\eps} F(\pj_t x)dW(t).  \label{eq:quadnon_main_in_short_form_diffinstab}
\end{align}
Writing  $\hprc(t):=\ham(\pj_tx)$ and $\hprc^\eps(t):=\hprc(t/\eps)$, we can show that the distribution of $\hprc^\eps$ converges weakly to the distribution of $$d\hlev_t=\Psiz G_q(\Phi \hlev_t)dt+|\Psiz F(\Phi \hlev_t)|dW_t.$$ However, for practical use, one might want to approximate $\hprc^\eps$ for small $\eps$ with $\hlev$. In this case, the following equation \emph{might} give a better approximation.
$$d\hlev_t=\Psiz G_q(\Phi\hlev_t)dt+ \eps b_H^{q,(2)}(\hlev_t) dt +\eps \Psiz G(\Phi\hlev_t)dt + |\Psiz F(\Phi\hlev_t)|dW_t,$$
where $b_H^{q,(2)}$ is given in \eqref{eq:HuttBifScenario_whiteNoise_cent}.
\end{rmk}

\cite{lefebvre} considers scalar systems satisfying assumption \ref{ass:assumptondetsys_Hutt}, but does not impose \eqref{eq:assumption_on_Gq_zero_Hutt}. \cite{lefebvre} gives a method to construct higher order corrections to the center-manifold in presence of periodic forcing and white noise. They show that having higher order corrections in the center-manifold would improve accuracy of reconstructing the trajectories (figures 2 and 6 in \cite{lefebvre}). However, these corrections should be evaluated through numerical simulations of a \emph{delay equation}---for example, the correction to the center-manifold in equation 52 of \cite{lefebvre} should be numerically simulated. In scalar equations this task can be circumvented by employing series solutions as in equation 53 of \cite{lefebvre}. However, for multidimensional system this involves evaluating reasonable number of eigenvalues and eigenvectors of the linear delay system. Further, the computations require memory for storing the history of Brownian motion for computing the convolutions (equation 55 in \cite{lefebvre}). The extra effort required from the methods in \cite{lefebvre} allows to reconstruct trajectories. The averaging methods presented in our article would deal with \emph{distributions alone} in the limit of small $\eps$ and cannot reconstruct trajectories.

Finally, for completeness, we consider equations of the form \eqref{eq:detDDE_pert_gennoise_nonon} with assumption \ref{ass:assumptondetsys_Hutt}. In this case it can be shown that theorem \ref{thm:gennoiseweakconvg} holds with 
\begin{align}\label{eq:HuttBifScenario_GenNoise_DriDiff}
b_H(h)=\left(\int_0^\infty R(s)ds\right)\left(\Ind F(\Phi h).\nabla\right)\Psiz F(\Phi h), \quad \,\, \sigma_H^2(h)=2\left(\int_0^\infty R(s)ds\right)\left(\Psiz F(\Phi h)\right)^2.
\end{align}


\section*{Acknowledgements}
The authors would like to gratefully acknowledge the suggestions of Prof. Volker Wihstutz. 

The work was financially supported by the National Science Foundation under grant numbers CMMI 1000906 and 1030144. Any opinions, findings, and
conclusions or recommendations expressed in this paper are those of the authors
and do not necessarily reflect the views of the National Science Foundation.




\appendix


\section{\label{appsec:PREpap:errorsofKuskeGaudFof} Errors in \cite{Klosek_Kuske_multi_siam}, \cite{Kuske_SD}, \cite{PhysRevE.85.056214} and shortcomings in \cite{Fofan02_sampstabchatter}, \cite{Fofana_PEM_chatter}, \cite{Caihong}.}


\subsection{\label{appsec:PREpap:errorsofKuskeGaudFof__Kusk} Errors in \cite{Klosek_Kuske_multi_siam}, \cite{Kuske_SD}}

One of the equations considered in \cite{Kuske_SD} is:
\begin{align}\label{eq:Kusk:main}
d\Xeps(t)=\frac{1}{\eps^2}\bigg(-\alpha \Xeps(t)+\beta \Xeps(t-\eps^2\tau)\bigg)dt+\Xeps(t)dW(t),
\end{align}
where $W$ is a Wiener process\footnote{This is time-rescaled version of eq 1.1 in \cite{Kuske_SD}. The analysis below appears in section 2 of \cite{Kuske_SD}.}.
The above system is studied as a perturbation of the linear system 
\begin{align}\label{eq:Kusk:main:unpert}
\dot{x}(t)=\frac{1}{\eps^2}\bigg(-\alpha x(t)+\beta x(t-\eps^2\tau)\bigg).
\end{align}
Seeking solution of the form $e^{\lambda t/\eps^2}$ the  characteristic equation is found to be $\lambda=-\alpha+\beta e^{-\lambda \tau}$. Let the parameters $\alpha,\beta, \tau=\tau_c+\eps^2\tau_2$ be such that when $\tau_2=0$, a pair of roots $\pm i\om$ are on the imaginary axis and all other roots are with negative real part. In this scenario we have $i\om = -\alpha + \beta e^{-i\om \tau_c}$ which on solving gives\footnote{This is eq 2.1 in \cite{Kuske_SD}.}
\begin{align}\label{eq:Kusk:paramatcrit}
\om=\sqrt{\beta^2-\alpha^2}, \qquad \beta \cos(\om \tau_c)=\alpha, \qquad \beta \sin(\om \tau_c)=-\om.
\end{align}

\cite{Kuske_SD} employs multiscale analysis and for that purpose writes\footnote{This is eq 2.11 in \cite{Kuske_SD}.}
\begin{align}\label{eq:Kusk:Wsplit}
dW(t)=\mkK_0dW_0(t)+\mkK_{2,1}\cos(\frac{2\om t}{\eps^2})dW_{2,1}(t)+\mkK_{2,2}\sin(\frac{2\om t}{\eps^2})dW_{2,2}(t),
\end{align}
where $W_{i}$ are independent Brownian motions.
\cite{Kuske_SD} assumes that solution $\Xeps$ is of the form\footnote{This is eq 2.2 in \cite{Kuske_SD}.}
\begin{align}\label{eq:Kusk:assumed_form}
\Xeps(t)=A(t)\cos(\om t/\eps^2)+B(t)\sin(\om t/\eps^2).
\end{align}
Here $A,B$ vary at different scale (in the spirit of multiscale analysis) than cosine and sine.

According to \cite{Kuske_SD}, on one hand, applying Ito formula we have\footnote{This is eq 2.4 in \cite{Kuske_SD}.}
\begin{align}\label{eq:Kusk:assumed_form_ito}
d\Xeps= \frac{1}{\eps^2}\left(-\om \mfs A + \om \mfc B\right)dt+ \mfc dA+ \mfs dB,
\end{align}
where $\mfc=\cos(\om t/\eps^2)$ and $\mfs=\sin(\om t/\eps^2)$.
On the other hand, since $\Xeps$ must satisfy \eqref{eq:Kusk:main} we must have\footnote{This is eq 2.5 in \cite{Kuske_SD}.}
\begin{align}\label{eq:Kusk:assumed_form_prop}
d\Xeps &= \frac{1}{\eps^2}\left(-\alpha\left(\mfc A+\mfs B\right)+\beta\left(A_\tau \cos(\frac{\om(t-\eps^2\tau)}{\eps^2})+B_\tau \sin(\frac{\om(t-\eps^2\tau)}{\eps^2})\right)\right)dt  \notag \\
& \qquad + (\mfc A+\mfs B)(\mkK_0dW_0(t)+\mkK_{2,1}\cos(\frac{2\om t}{\eps^2})dW_{2,1}(t)+\mkK_{2,2}\sin(\frac{2\om t}{\eps^2})dW_{2,2}(t)),
\end{align}
where $A_\tau$ means $A(t-\eps^2\tau)$. 

Using $\tau=\tau_c+\eps^2\tau_2$ and \eqref{eq:Kusk:paramatcrit} we have
\begin{align}\label{eq:Kusk:cossinapprx}
\beta \cos(\frac{\om(t-\eps^2\tau)}{\eps^2})& = (\alpha \mfc -\om \mfs)+ \eps^2 \om \tau_2(\om \mfc + \alpha \mfs) \\
\beta \sin(\frac{\om(t-\eps^2\tau)}{\eps^2})& = (\om \mfc +\alpha \mfs)+ \eps^2 \om \tau_2(-\alpha \mfc + \om \mfs). 
\end{align}
Using the above in \eqref{eq:Kusk:assumed_form_prop} and comparing the resulting equation with \eqref{eq:Kusk:assumed_form_ito} we have
\begin{align}\label{eq:Kusk:ABnoap_pre}
&\frac{1}{\eps^2}\left(-\alpha(\mfc A+\mfs B)+A_\tau(\alpha \mfc-\om \mfs)+B_\tau(\alpha \mfs + \om \mfc )\right)dt \\
& \qquad + \om \tau_2\left(\om(\mfc A_\tau+\mfs B_\tau)+\alpha(\mfs A_\tau-\mfc B_\tau)\right)dt \notag \\
& \qquad + (\mfc A+\mfs B)\left(\mkK_0dW_0(t)+\mkK_{2,1}\cos(\frac{2\om t}{\eps^2})dW_0(t)+\mkK_{2,2}\sin(\frac{2\om t}{\eps^2})dW_0(t)\right) \notag \\
& \qquad -\frac{1}{\eps^2}\left(-\om \mfs A + \om \mfc B\right)dt- \mfc dA- \mfs dB \qquad =\qquad 0. \notag 
\end{align}
\cite{Kuske_SD} then multiplies the above with $\mfc$ or $\mfs$ and integrates over a time period, while treating $A$ and $B$ as constants, to get the following equations:
\begin{align}\label{eq:Kusk:ABnoap}
dA&=-\alpha \hat{d}A - \om \hat{d}B+\om \tau_2(\om A_\tau - \alpha B_\tau)dt + A\mkK_{2,0}dW_0+ \frac12A\mkK_{2,1}dW_{2,1}+ \frac12B\mkK_{2,2}dW_2 \notag \\
dB&=\om \hat{d}A - \alpha \hat{d}B+\om \tau_2(\alpha A_\tau + \om B_\tau)dt + B\mkK_{2,0}dW_0- \frac12B\mkK_{2,1}dW_{2,1}+ \frac12A\mkK_{2,2}dW_2,
\end{align}
where $\hat{d}A$ means $\frac{A(t)-A(t-\eps^2\tau)}{\eps^2}dt$.

In \eqref{eq:Kusk:ABnoap} the constants $\mkK$ are not yet determined.  \cite{Kuske_SD} determines them in the following way: \cite{Kuske_SD} compares the diffusive part of the generator for $\Xeps$ and for $(A,B)$. The diffusive part of the generator for $(A,B)$ is
\begin{align}\label{eq:Kusk:diffusivepartofgenforAB}
&(A^2\partial_A\partial_A+B^2\partial_B\partial_B+2AB\partial_A\partial_B)\mkK_{2,0}^2 \notag\\
&\qquad + \frac{1}{4}(A^2\partial_A\partial_A+B^2\partial_B\partial_B-2AB\partial_A\partial_B)\mkK_{2,1}^2 \notag \\
&\qquad + \frac{1}{4}(B^2\partial_A\partial_A+A^2\partial_B\partial_B+2AB\partial_A\partial_B)\mkK_{2,2}^2.
\end{align}
The diffusive part of the generator for $x$ is
\begin{align}\label{eq:Kusk:diffusivepartofgenforx}
x^2\partial_x\partial_x = (\mfc A+\mfs B)^2(\mfc \partial_A+\mfs \partial_B)^2.
\end{align}
Averaging \eqref{eq:Kusk:diffusivepartofgenforx} over one time period, \cite{Kuske_SD} obtains\footnote{This is eq 2.16 in \cite{Kuske_SD}.}
\begin{align}\label{eq:Kusk:diffusivepartofgenforxAvg}
\frac{3A^2+B^2}{8}\partial_A\partial_A+ \frac{3B^2+A^2}{8}\partial_B\partial_B+\frac12AB\partial_A\partial_B.
\end{align}
\cite{Kuske_SD} equates \eqref{eq:Kusk:diffusivepartofgenforxAvg} and \eqref{eq:Kusk:diffusivepartofgenforAB} to find that 
\begin{align}\label{eq:Kusk:equategens}
\mkK_{2,0}=\frac12, \qquad \mkK_{2,1}=\mkK_{2,2}=\frac{1}{\sqrt{2}}.
\end{align}
Then \cite{Kuske_SD} presents a figure showing that density of $A(T)\cos(\om T/\eps^2)+B(T)\sin(\om T/\eps^2)$, with $A, B$ simulated from \eqref{eq:Kusk:ABnoap}, gives good approximation to the density of $\Xeps(T)$.

The above procedure is not convincing due to the following reasons:
\begin{itemize}
\item It is not clear whether the error in transferring from \eqref{eq:Kusk:ABnoap_pre} to \eqref{eq:Kusk:ABnoap} would go to zero in some sense as $\eps \to 0$. 
\item Note that \eqref{eq:Kusk:ABnoap} is still a delay equation and hence there would not be much advantage in simulating $A,B$ compared to simulating $\Xeps$. The delay itself is small $O(\eps^2)$, but the difference $A(t)-A(t-\eps^2\tau)$ is magnified by $\eps^{-2}$.
\item Note that, heuristically, the LHS of \eqref{eq:Kusk:Wsplit} is a normal random variable with variance $dt$; and hence, for consistency, we must have
\begin{align}\label{eq:Kusk:Wsplit_consistency}
\mkK_{2,0}^2+\mkK_{2,1}^2\cos^2(\frac{2\om t}{\eps^2})+\mkK_{2,2}^2\sin^2(\frac{2\om t}{\eps^2})=1.
\end{align}
The above is possible only if we take $|\mkK_{2,1}|=|\mkK_{2,2}|$ and set 
\begin{align}\label{eq:Kusk:Wsplit_consistency_2}
\mkK_{2,0}^2+\mkK_{2,1}^2=1.
\end{align}
But note that {\bf{\eqref{eq:Kusk:equategens} contradicts the consistency equation \eqref{eq:Kusk:Wsplit_consistency_2}}.} We have from \eqref{eq:Kusk:equategens} that $\mkK_{2,0}^2+\mkK_{2,1}^2=\frac{3}{4}\neq 1$.
\end{itemize}

We show by means of numerical simulation that the above procedure is indeed wrong.

 In \eqref{eq:Kusk:main} set $\alpha=0$, $\beta=-\frac{\pi}{2}$ and $\tau_c=1$, $\tau_2=0$. Then $\om=\frac{\pi}{2}$ and this system satisfies assumption \ref{ass:assumptondetsys}. The equations \eqref{eq:Kusk:ABnoap} in this case becomes:
\begin{align}\label{eq:Kusk:ABnoap_specialex}
\left(\begin{array}{c}dA\\dB\end{array}\right)&=\frac{1}{\eps^2}\left(\begin{array}{cc}0 & -\om \\ \om & 0\end{array}\right)\left(\begin{array}{c}A(t)-A(t-\eps^2)\\B(t)-B(t-\eps^2)\end{array}\right)dt\\ \notag
& \qquad + \frac12\left(\begin{array}{cc}1 & 0 \\ 0 & 1\end{array}\right)\left(\begin{array}{c}A(t)\\B(t)\end{array}\right)dW_{2,0}  \\ \notag
& \qquad + \frac{1}{2\sqrt{2}}\left(\begin{array}{cc}1 & 0 \\ 0 & -1\end{array}\right)\left(\begin{array}{c}A(t)\\B(t)\end{array}\right)dW_{2,1} \\ \notag
& \qquad +  \frac{1}{2\sqrt{2}}\left(\begin{array}{cc}0 & 1 \\ 1 & 0\end{array}\right)\left(\begin{array}{c}A(t)\\B(t)\end{array}\right)dW_{2,2}
\end{align}

Numerical simulations show that splitting $W$ into harmonics as in \eqref{eq:Kusk:Wsplit} is unnecessary. For this purpose, consider 
\begin{align}\label{eq:Kusk:ABnoap_specialex_corrected}
\left(\begin{array}{c}dA\\dB\end{array}\right)&=\frac{1}{\eps^2}\left(\begin{array}{cc}0 & -\om \\ \om & 0\end{array}\right)\left(\begin{array}{c}A(t)-A(t-\eps^2)\\B(t)-B(t-\eps^2)\end{array}\right)dt\\ \notag
& \qquad + \left(\begin{array}{cc}1 & 0 \\ 0 & 1\end{array}\right)\left(\begin{array}{c}A(t)\\B(t)\end{array}\right)dW_{2,0}.
\end{align}
i.e.  $\mkK_{0}=1$, $\mkK_{2,1}=0=\mkK_{2,2}$.

\begin{center}
\begin{figure}
\centering
\includegraphics[scale=0.7]{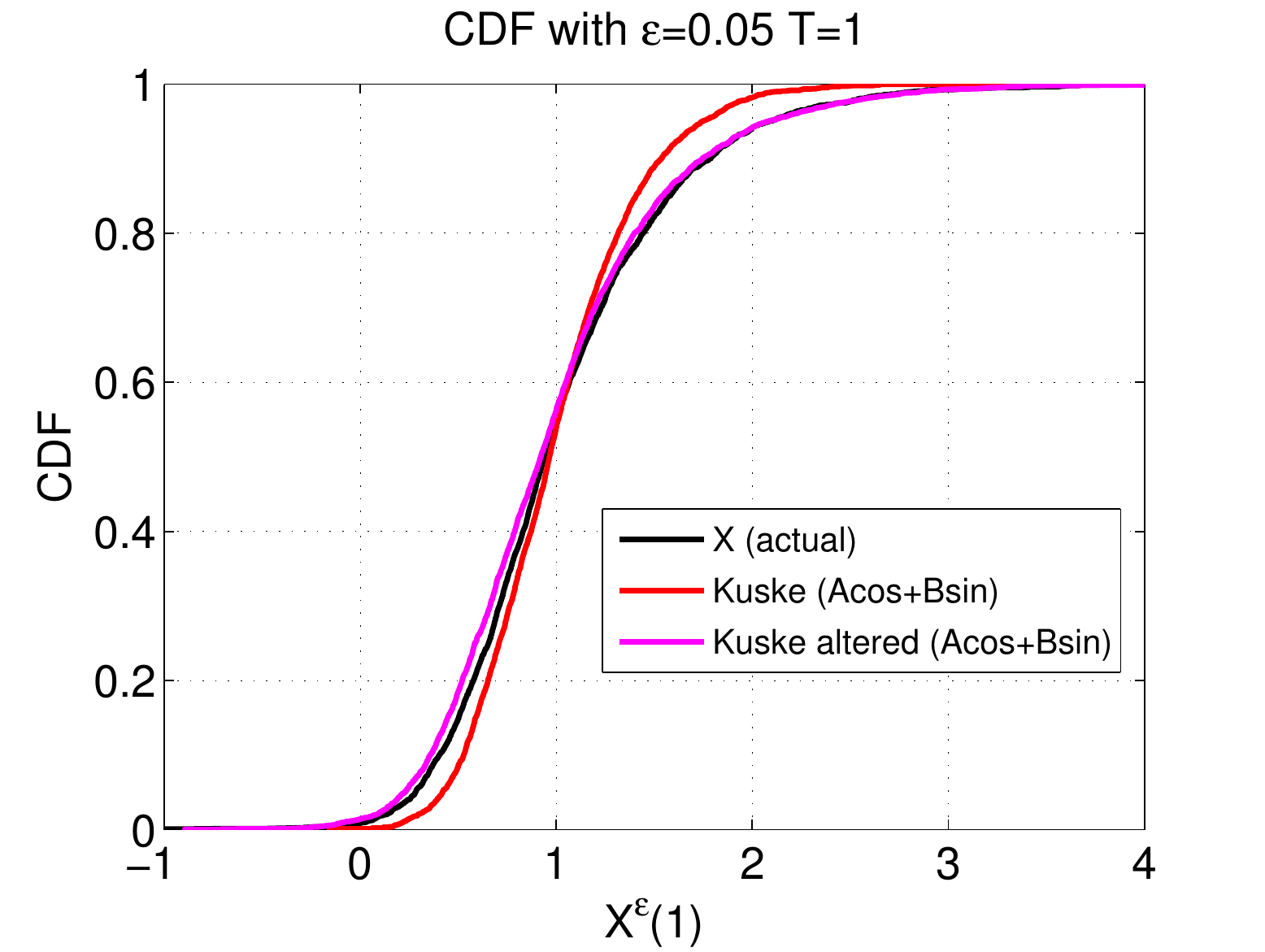}
\caption{\emph{X (actual)} is obtained from simulating the original dynamics  \eqref{eq:Kusk:main}. \emph{Kuske (Acos+Bsin)} is $A(T)\cos(\om T/\eps^2)+B\sin(\om T/\eps^2)$ obtained from simulating \eqref{eq:Kusk:ABnoap_specialex}. \emph{Kuske altered (Acos+Bsin)} is $A(T)\cos(\om T/\eps^2)+B\sin(\om T/\eps^2)$ obtained from simulating \eqref{eq:Kusk:ABnoap_specialex_corrected}.}
\label{fig:KuskWrg1X}
\end{figure}
\end{center}

We set $\eps=0.05$, $T=1$. The initial condition is $\Xeps(t)=\cos(\om t/\eps^2)$ for $t\in [-\eps^2,0]$, i.e. $\pjeps_0\Xeps(\theta)=\cos(\om \theta)$ for $\theta \in [-1,0]$, i.e. $A(t)=1$ for $t\leq 0$ and $B(t)=0$ for $t\leq 0$. The cumulative distribution in the figure \ref{fig:KuskWrg1X} is obtained with 2400 realizations.

Figure \ref{fig:KuskWrg1X} shows that \eqref{eq:Kusk:ABnoap_specialex_corrected} better matches the actual dynamics \eqref{eq:Kusk:main} than \eqref{eq:Kusk:ABnoap_specialex}. But, note that \eqref{eq:Kusk:ABnoap_specialex_corrected} is still a delay equation and there is no advantage in simulating $(A,B)$ compared to simulating $X$. 

\subsection{\label{appsec:PREpap:errorsofKuskeGaudFof__Gaud} Errors in \cite{PhysRevE.85.056214} and \cite{Gaud_hopf_DDE_multi_scale}}

There are two errors in the analysis of \cite{PhysRevE.85.056214} and \cite{Gaud_hopf_DDE_multi_scale}, one of which is similar in nature to the previous section. We illustrate the errors using a special case of the equation considered in \cite{PhysRevE.85.056214}.

\cite{PhysRevE.85.056214} considers
\begin{align}\label{eq:Gaud:main_sp}
\ddot{x}(t)+x(t)+\eta x(t-1)-\beta \dot{x}(t)=\sqrt{2D} x(t)\xi(t),
\end{align}
where $\xi$ is a white noise process with correlation $\expt[\xi(t)\xi(t')]=\delta(t-t')$. For now, lets set $D=0$. The characterisitc equation is $\lambda^2+1+\eta e^{-\lambda}-\beta\lambda=0$. Given $\eta$, solve $\eta\cos\om = \om^2-1$ for $\om$ and get $\beta_c=-\eta \sin\om/\om$. With $\beta=\beta_c$ the system \eqref{eq:Gaud:main_sp} (with $D=0$) satisfies assumption \ref{ass:assumptondetsys} with critical roots of the characteristic equation being $\pm i\om$. We assume $\beta=\beta_c$.

\cite{PhysRevE.85.056214} assumes the solution is of the form
\begin{align}\label{eq:Gaud:solnform_ms}
x(t,T)=\eps A(T)\cos\om t\,-\,\eps B(T)\sin\om t
\end{align}
where $T=\eps^2t$ is the slow time scale. Then, 
\begin{align}\label{eq:Gaud:solnform_ms_delayterms}
x(t-1,T-\eps^2)&=x(t,T)\cos\om - (\sin\om/\om) \partial_tx(t,T) \\ \notag
& \qquad \qquad -\eps^2\eps\frac{A(T)-A(T-\eps^2)}{\eps^2}\cos(\om(t-1)) \\ \notag
& \qquad \qquad +\eps^2\eps\frac{B(T)-B(T-\eps^2)}{\eps^2}\sin(\om(t-1)).
\end{align}
But, \cite{PhysRevE.85.056214} sets last two terms in the RHS to zero claiming $A(T)\approx A(T-\eps^2)$ and $B(T)\approx B(T-\eps^2)$. 
{\bf However, as $\eps \to 0$ it is easy to see that (if derivative of $A$ and $B$ exist) these terms go to $\partial_TA$ and $\partial_TB$ respectively. At which $\eps$ should we ignore these and which $\eps$ should we consider it as a derivative?}

Differentiating, we get
\begin{align}\label{eq:Gaud:solnform_ms_xdot}
\dot{x}(t)& = (\eps^2\partial_T+\partial_t)x(t,T)=\eps^2(\eps\partial_TA\cos\om t-\eps \partial_TB\sin\om t) + \partial_tx(t,T)
\end{align}
\begin{align}\label{eq:Gaud:solnform_ms_xddot}
\ddot{x}(t)&= (\eps^2\partial_T+\partial_t)^2x(t,T) =\eps^4(\eps\partial_T^2A\cos\om t-\eps \partial_T^2B\sin\om t)\\ \notag
& \qquad -\eps^22\om(\eps\partial_TA\sin\om t + \eps \partial_TB\cos\om t)-\om^2x(t,T)
\end{align}

Putting \eqref{eq:Gaud:solnform_ms_delayterms}, \eqref{eq:Gaud:solnform_ms_xdot} and \eqref{eq:Gaud:solnform_ms_xddot} together in \eqref{eq:Gaud:main_sp} and using $\eta\cos\om = \om^2-1$,  $\beta_c=-\eta \sin\om/\om$ and ignoring terms of order more than $\eps^3$ we get that
\begin{align}\label{eq:Gaud:postmswhatactual}
-2\om\eps^3 & (\partial_TA\sin\om t  +\partial_TB\cos\om t)\\ \notag
& -\eps^3\eta(\Delta A(T)\cos(\om(t-1))-\Delta B(T)\sin(\om(t-1))) \\ \notag
& -\eps^3\beta_c(\partial_TA\cos\om t-\partial_TB\sin\om t) \qquad = \\ \notag
& \qquad \qquad = \sqrt{2D} \eps \bigg(A(T)\cos\om t-B(T)\sin\om t\bigg)\xi(t),
\end{align}
where $\Delta A(T)$ means $\frac{A(T)-A(T-\eps^2)}{\eps^2}$ etc.
The corresponding equation that \cite{PhysRevE.85.056214} arrives at\footnote{This is equation 9 in \cite{PhysRevE.85.056214}. The quantity $\mu$ defined under equation 7 of \cite{PhysRevE.85.056214} is zero for the special case that we consider.} is:
\begin{align}\label{eq:Gaud:postmswhatGaudgets}
-\om\eps^3 & (\partial_TA\sin\om t  +\partial_TB\cos\om t)\\  \notag
& \qquad \qquad = \sqrt{2D} \eps \bigg(A(T)\cos\om t-B(T)\sin\om t\bigg)\xi(t),
\end{align}
The equation \eqref{eq:Gaud:postmswhatGaudgets} does not match with \eqref{eq:Gaud:postmswhatactual} when $\Delta A$, $\Delta B$ are set to zero, nor when they are set as actual derivatives $\partial_TA$, $\partial_TB$.

\cite{PhysRevE.85.056214} proceeds with \eqref{eq:Gaud:postmswhatGaudgets}, multiplies with $\sin\om t$ and averages over a time period to arrive at:
\begin{align}\label{eq:Gaud:postmswhatGaudgets_preavg}
-\om\eps^3 \frac12 \partial_TA &= \sqrt{2D} \eps \bigg(A(T)\llb \cos\om t \sin\om t \,\xi(t)\rrb-B(T)\llb \sin^2\om t\, \xi(t)\rrb\bigg), \\ \notag
&= \sqrt{2D} \eps \frac12\bigg(A(T)\llb \sin2\om t \,\xi(t)\rrb-B(T)\llb \xi(t)\rrb+B(T)\llb \cos2\om t \,\xi(t)\rrb\bigg),
\end{align}
where $\llb \,\rrb$ is used for \emph{time-averaging}.

The intermediate steps in \cite{PhysRevE.85.056214} are not clear, but the end result of \cite{PhysRevE.85.056214} is that $D$ is scaled as $D=\eps^2\tilde{D}$ and three new Gaussian process $\xi_0,\xi_1,\xi_2$ are defined on slow time scale and the following are used:
\begin{align}\label{eq:Gaud:newGauss}
\llb \xi(t)\rrb =\eps \xi_0, \qquad \llb \cos2\om t  \,\xi(t)\rrb = \frac{\eps}{\sqrt{2}}\xi_1, \quad \llb \sin2\om t  \,\xi(t)\rrb = \frac{\eps}{\sqrt{2}}\xi_2.
\end{align}
Employing this in \eqref{eq:Gaud:postmswhatGaudgets_preavg} the following is arrived at:
\begin{align}\label{eq:Gaud:eqforA}
-\frac{\om}{\sqrt{2\tilde{D}}} \partial_TA = -B\xi_0+\frac{1}{\sqrt{2}}B\xi_1+\frac{1}{\sqrt{2}}A\xi_2.
\end{align}
Similary, \cite{PhysRevE.85.056214} multiplies \eqref{eq:Gaud:postmswhatGaudgets} with $\cos\om t$ and averages over a time period and employs \eqref{eq:Gaud:newGauss} to arrive at:
\begin{align}\label{eq:Gaud:eqforB}
-\frac{\om}{\sqrt{2\tilde{D}}}  \partial_TB = A\xi_0+\frac{1}{\sqrt{2}}A\xi_1-\frac{1}{\sqrt{2}}B\xi_2.
\end{align}
The equations \eqref{eq:Gaud:eqforA} and \eqref{eq:Gaud:eqforB} are respectively (16) and (17) in \cite{PhysRevE.85.056214}.

Now we show that the above method is not consistent with itself. From \eqref{eq:Gaud:eqforA} and \eqref{eq:Gaud:eqforB} we get
\begin{align}\label{eq:Gaud:postmswhatGaudgets_inconsit}
-\frac{\om}{\sqrt{2\tilde{D}}} & (\partial_TA\sin\om t  +\partial_TB\cos\om t)\\  \notag
& \qquad \qquad = (-B\mfs+A\mfc)\xi_0+\frac{1}{\sqrt{2}}(B\mfs+A\mfc)\xi_1+\frac{1}{\sqrt{2}}(A\mfs-B\mfc)\xi_2,\\
& \qquad \qquad =: \mathfrak{F}(T)
\end{align}
where $\mfs=\sin\om t$ and $\mfc=\cos\om t$. Now $\expt[ \mathfrak{F}(T) \mathfrak{F}(T)]$ equals
\begin{align}\label{eq:Gaud:postmswhatGaudgets_varcalc}
(-B\mfs+A\mfc)^2+\frac12(B\mfs+A\mfc)^2 & +\frac12(A\mfs-B\mfc)^2 \\ \notag
& =(A\mfc-B\mfs)^2+\frac12(A^2+B^2).
\end{align}
But from \eqref{eq:Gaud:postmswhatGaudgets}
\begin{align}\label{eq:Gaud:postmswhatGaudgets_comparevar}
-\frac{\om}{\sqrt{2\tilde{D}}}&(\partial_TA\sin\om t  +\partial_TB\cos\om t)\\  \notag
& \qquad \qquad =  \eps \big(A\mfc-B\mfs\big)\xi(t) =: \eps \mathfrak{F}(T),
\end{align}
Now $\expt[ \mathfrak{F}(T) \mathfrak{F}(T)]$ equals $(A\mfc-B\mfs)^2$. So the system \eqref{eq:Gaud:eqforA},\eqref{eq:Gaud:eqforB} has an extra variance of $\frac12(A^2+B^2)$ (see \eqref{eq:Gaud:postmswhatGaudgets_varcalc}) than what is required.

\subsection{\label{appsec:PREpap:errorsofKuskeGaudFof__Fof} Shortcomings in \cite{Fofan02_sampstabchatter}, \cite{Fofana_PEM_chatter}, \cite{Caihong}}

\cite{Fofan02_sampstabchatter}, \cite{Fofana_PEM_chatter} consider oscillators that arise in machine tool dynamics and \cite{Caihong} considers human standing model. They apply the spectral theory of linear DDE just like is done in this paper. However, right from the beginning of the analysis they claim that the stable ($Q$) part of the solution can be ignored. They take noise as Wiener process but do not consider stronger deterministic perturbations $G_q$ as in section \ref{sec:strongerpert}. However, when considering $G_q$ or when considering other noise processes, ignoring the $Q$ part of the solution would lead to wrong results. As pointed out in remarks \ref{rmk:Foferror1} and \ref{rmk:Foferror2}, this leads to loss of some of the drift terms.


\section{\label{appsec:strongdetpert_example_nodel} An example illustrating the approach for calculation of $b_H^{q,(i)}$ in theorem \ref{thm:PREpap:convgHforWquad}}
Consider the system without delay given by $\ddot{x}+x=\eps \dot{x}y$, and $\dot{y}=-y+\eps\dot{x}^2$. Here $x$ is oscillatory and $y$ is stable. The quantity $\hprc=\frac12(x^2+\dot{x}^2)$ evolves slowly compared to $x$ and $y$. Writing in state-space form $z_1=x$, $z_2=\dot{x}$ we have
\begin{align}\label{appsec:strongdetpert_example_nodel:example_eq_statspace}
\left(\begin{array}{c}\dot{z_1}\\ \dot{z_2} \\ \dot{y}\end{array}\right)=\left(\begin{array}{c}z_2\\ -z_1 \\ -y\end{array}\right)+\eps\left(\begin{array}{c}0 \\ z_2y \\ z_2^2\end{array}\right)
\end{align}
and $\dot{\hprc}=\eps\hprcdrift^{(q)}(z,y)$, where $\hprcdrift^{(q)}(z,y)=z_2^2y$. 

The unperturbed system is obtained by setting $\eps=0$ in \eqref{appsec:strongdetpert_example_nodel:example_eq_statspace}. The differential of any function $f$ along trajectory of unperturbed system is given by $\igen_0f$  where $\igen_0=z_2\frac{\partial}{\partial z_1}-z_1\frac{\partial}{\partial z_2}-y\frac{\partial}{\partial y}$. The differential along the perturbations is given by $\igen_1f$  where $\igen_1=z_2y\frac{\partial}{\partial z_2}+z_2^2\frac{\partial}{\partial y}$. Note that $\dot{f}(z_t,y_t)=((\igen_0+\eps\igen_1)f)(z_t,y_t)$.

Now let 
\begin{align}\label{eq:quadcorr_ptf}
H(z,y)=\hprc(z)-\eps c(z,y)+\eps^2g_1(z,y)+\eps^2g_2(z)
\end{align}
 where $c,g$ are yet to be determined. On differentiating we get (until order $\eps^2$)
\begin{align}\label{appsec:strongdetpert_example_nodel:example_eq_Hexpansion}
\dot{H}(z_t,y_t)=\eps\big(\hprcdrift^{(q)}(z_t,y_t)-(\igen_0c)(z_t,y_t)\big)-\eps^2(\igen_1c)(z_t,y_t)+\eps^2(\igen_0g_1)(z_t,y_t)+\eps^2(\igen_0g_2)(z_t,y_t) + O(\eps^3).
\end{align}
Now, choose $c$ such that $\igen_0c=\hprcdrift^{(q)}$. Choose $g_1$ such that $(\igen_0g_1)(z,y)=(\igen_1c)(z,y)-(\igen_1c)(z,0)$. Such a choice of $g_1$ is possible because, according to the unperturbed dynamics $y$ decays to zero exponentially fast. Now, note that $(\igen_1c)(z,0)$ is a function of $z$ alone; and  the unperturbed $z$ dynamics is `oscillation with constant amplitude $\sqrt{2\hprc}$'.
Now, let the average of $(\igen_1c)(z,0)$ along an orbit of constant $\hprc$ be denoted by $\{\igen_1 c\}$. This $\{\igen_1 c\}$ would be a function only of $\frac12(z_1^2+z_2^2)$ or what is the same --- $\hprc$.
Choose $g_2(z)$ such that $(\igen_0g_2)(z,0)=(\igen_1c)(z,0)-\{\igen_1 c\}|_{\frac12(z_1^2+z_2^2)}$. Plugging the above choices of functions in \eqref{appsec:strongdetpert_example_nodel:example_eq_Hexpansion} we get
\begin{align}\label{appsec:strongdetpert_example_nodel:example_eq_HexpansionFin}
\dot{H}(z_t,y_t)=-\eps^2\{\igen_1 c\}|_{\hprc} + O(\eps^3).
\end{align}
Hence, for times of order $O(1/\eps^2)$ we have $H(z_t,y_t)=H(z_0,y_0)+\eps^2\int_0^t\{\igen_1 c\}|_{\hprc_s}ds+O(\eps)$. Since $H$ differs from $\hprc$ only by $O(\eps)$ (see \eqref{eq:quadcorr_ptf}) we can write $\hprc_t=\hprc_0+\eps^2\int_0^t\{\igen_1 c\}|_{\hprc_s}ds+O(\eps).$
So, for times of order $O(1/\eps^2)$, if we use 
\begin{align}\label{appsec:strongdetpert_example_nodel:example_eq_HexpansionFinEnd}
\dot{\hprc}=-\eps^2\{\igen_1 c\}|_{\hprc}
\end{align}
then the error resulted in $\hprc$ would be only of $O(\eps)$.
Such a method is shown in \cite{ChowParet}---we have adapted it to stochastic delay equations in \cite{LNNSNarxSD}. 

To see why the above method is useful, note that $c$ in $\igen_0c=\hprcdrift^{(q)}$ can be immediately solved using method of characterisitcs. Since the solution to the unperturbed system is $z_1(t)=z_1(0)\cos t + z_2(0)\sin t$, $z_2(t)=-z_1(0)\sin t + z_2(0)\cos t$, $y(t)=y(0)e^{-t}$, and $\hprcdrift^{(q)}(z,y)=z_2^2y$ we get
$c(z,y)=-\int_0^\infty(-z_1\sin t + z_2\cos t )^2ye^{-t}dt$. Now, $(\igen_1c)(z,0)=-\int_0^\infty z_2^2(-z_1\sin t + z_2\cos t )^2e^{-t}dt$. Hence $\{\igen_1 c\}|_{\hprc}$ is
\begin{align*}
\frac{1}{2\pi}\int_0^{2\pi}\left(-\int_0^\infty z_2^2(-z_1\sin t + z_2\cos t )^2e^{-t}dt\right)&\bigg|_{(z_1,z_2)=\sqrt{2\hprc}(\sin s, \cos s)}ds \\
&=-\hprc^2\int_0^\infty\frac12(2+\cos 2t)e^{-t}dt \quad & = -\frac{11}{10}\hprc^2.
\end{align*}
So we have $\dot{\hprc}=\eps^2\frac{11}{10}\hprc^2 + O(\eps^3)$. The reader can check using conventional center-manifold calculations that same answer would be obtained. However the method presented here would easily adapt to multidimensional delay equations as shown in section \ref{sec:strongerpert}.


\section{\label{subsec:PREpap:helpfulcomm} Explicit evaluation of $b_H^{q,(k)}$ using \eqref{eq:taudef}--\eqref{eq:avgbhqk}}

In this section we show how the explicit formulas \eqref{eqn:explictcentmanterms_form_1}--\eqref{eqn:explictcentmanterms_form_2} can be derived from \eqref{eq:taudef}--\eqref{eq:avgbhqk}. First we give a few preliminaries.

Recall that, for $\icond \in \C$, $\Th(t)\icond$ denotes the solution at time $t$ of the unperturbed linear system \eqref{eq:detDDE} with initial condition $\pj_0x=\icond$. Recall that $\C=P\oplus Q$ where $P$ is the space corresponding to the critical eigenvalues $\pm i\om_c$.
Recalling the evolution on $P$ defined by \eqref{eq:proj_eqns_unpert}, we have that for $\ubar{u}\in \mathbb{C}^2$ with $u_2=\bar{u_1}$, 
\begin{align}\label{eq:aux:solnonP}
\Th(t)\Phi \ubar{u}=\Phi e^{Bt}\ubar{u}.
\end{align}

Using \eqref{eq:aux:projofIndontoP} and \eqref{eq:aux:solnonP} we have for $n\times 1$ vector $\ubar{v}$
\begin{align}
\Th(t)\pi\Ind \ubar{v}=\Phi e^{Bt}\Psiz \ubar{v}.
\end{align}

For $\etavg_t$ defined in \eqref{eq:unpertsolused4avg}, we have $\Th(s)\etavg_t=\etavg_{t+s}$. The $z$ coordinates $\la \Psi,\Th(s)\etavg_{t} \ra$ are given by  $\frac12\sqrt{2\hlev} \left[\begin{array}{c}e^{i\om_c (t+s)} \\ e^{-i\om_c (t+s)}\end{array}\right]$ and hence for $\rho$ defined in \eqref{eq:taudef}, we can take $\rho(\etavg_t)=\frac{2\pi}{\om_c}-t$.

Using product rule for differentiation on $\hprcdrift^{q,(1)}$ (defined in \eqref{eq:hprcdrifdef_quad1})  and linearity of the function $E$, we have for $\xi,\eta \in \C$
\begin{align*}
(\xi.\nabla)\hprcdrift^{q,(1)}(\eta)&=E(\xi)G_q(\pi\eta) + E(\eta) (\pi\xi.\nabla)G_q(\pi\eta).
\end{align*}

Using product rule for differentiation on $\hprcdrift^{q,(2)}$ (defined in \eqref{eq:hprcdrifdef_quad2})  we have for $\xi,\eta \in \C$
\begin{align*}
(\xi.\nabla)\hprcdrift^{q,(2)}(\eta)&=\bigg((\xi.\nabla)E(\eta)\bigg)(G_q(\eta)-G_q(\pi\eta)) + E(\eta) (\xi.\nabla)G_q(\eta) - E(\eta)(\pi\xi.\nabla)G_q(\pi\eta).
\end{align*}
Since $\etavg_t$ (used in \eqref{eq:avgbhqk}) belongs to $P$, i.e. $\etavg_t=\pi\etavg_t$, the first term vanishes. Using linearity of differentials we have that
\begin{align}\label{eq:auxcalc4bhq2eval}
(\xi.\nabla)\hprcdrift^{q,(2)}(\etavg_t)&=E(\etavg_t)((I-\pi)\xi.\nabla)G_q(\etavg_t) \quad \text{ for all } \xi \in \C.
\end{align}

Now we show how \eqref{eqn:explictcentmanterms_form_2} can be derived. Using \eqref{eq:a2qdef_aux} in \eqref{eq:avgbhqk} we encounter the task of evaluating the differential $(\xi.\nabla)\hprcdrift^{q,(2)}(\Th(s)\etavg_t)$ with $\xi=\Th(s) \Ind G_q(\etavg_t)$. Using $\Th(s)\etavg_t=\etavg_{t+s}$ and \eqref{eq:auxcalc4bhq2eval} we get the differential as $E(\etavg_{t+s})((I-\pi)\Th(s)\Ind G_q(\etavg_t).\nabla)G_q(\etavg_{t+s})$. It is a property of the  unperturbed system that $\Th$ commutes with $(I-\pi)$. Defining $\mathcal{E}_t=e^{-i\om_ct}\Psiz_1+e^{i\om_ct}\Psiz_2$ we can write $E(\etavg_t)=\sqrt{2\hlev}\mathcal{E}_t$. So we can rewrite the differential as $\sqrt{2\hlev}(\Th(s)(I-\pi)\Ind G_q(\etavg_t).\nabla)(\mathcal{E}_{t+s}G_q(\etavg_{t+s}))$. Writing $G_q(\etavg_t)=\sum_{j=1}^n(G_q(\etavg_t))_j\ubar{e_j}$ and using linearity of differentials we get the desired form in \eqref{eqn:explictcentmanterms_form_2}.

\eqref{eqn:explictcentmanterms_form_1} can be similarly derived.


\section{\label{appsec:sketchproofGenNoise} A sketch of proof of theorem \ref{thm:gennoiseweakconvg}}
One way to characterize the probability distribution of a stochastic process $Y$ is by an operator called the infinitesimal generator $\igen$ defined as follows: for any nice real-valued function $f$ of the process $Y$, 
\begin{align}\label{eq:def:infgendef_illus}
(\igen f)(y)\overset{\text{def}}=\lim_{t \to 0}\frac1t({\expt[f(Y_{t})|Y_0=y]-f(y)}).
\end{align}
Here the `$\expt$' term means ``the average of $f(Y_t)$ given that the initial condition $Y_0$ equals $y$''. 
For example, the process whose infinitesimal generator is defined by $\igen f=\frac12f''$ has the same probability distribution as the standard Brownian motion. The process whose infinitesimal generator is defined by $(\igen f)(y)=b(y)f'(y)+\frac12\sigma^2(y) f''(y)$ has the same probability distribution as the process governed by the SDE, $dY=b(Y)dt+\sigma(Y)dW$ with $W$ a Wiener process. The process whose infinitesimal generator is $(\igen f)(y)=b(y)f'(y)$ is the ordinary differential equation $\dot{Y}=b(Y)$. \emph{The infinitesimal generator characterizes the probability distribution of a process.}

We consider the system \eqref{eq:detDDE_pert_rescale_gennoise}--\eqref{eq:evolofhprceps_gennoise} and try to find the infinitesimal generator $\igen_H$ of the process $\lim_{\eps \to 0}\hprc^{\eps}$. For this purpose consider the triplet process $(\pjeps X,\gnoiseeps,\hprc^\eps)$. It has the infinitesimal generator $\igen^\eps=\frac{1}{\eps^2}\igen_0+\frac{1}{\eps}\igen_1$, where for function $f$ of $(\eta,\gnoise,h)$
\begin{align}\label{eq:infgen_gennoise_L0_illus}
(\igen_0f)(\eta,\gnoise,h)&=(\noisegent f)(\eta,\gnoise,h) + \frac{d}{dt}\big|_{t=0}f(\Th(t)\eta,\gnoise,h),\\ \label{eq:infgen_gennoise_L1_illus}
(\igen_1f)(\eta,\gnoise,h)&=\sigma(\gnoise)(\Ind F(\eta).\nabla)f(\eta,\gnoise,h) + \sigma(\gnoise)\hprcdrift(\eta)\frac{\partial f}{\partial h}(\eta,\gnoise,h).
\end{align}
Here $\noisegent$ is the infinitesimal generator of the noise process $\gnoise$. Recall that $\Th(t)\eta$ is the solution at time $t$ of the unperturbed system \eqref{eq:detDDE} with initial condition $\eta$, and $\Ind$ is the matrix valued function defined in \eqref{eq:matrixInd0}. 

The following comments help in gaining an insight into the structure of $\igen^\eps$. Consider \eqref{eq:detDDE_pert_rescale_gennoise}--\eqref{eq:evolofhprceps_gennoise}. If there were no noise perturbations at all, then $\hprc^\eps$ would have remained a constant and $\pjeps \Xeps$ would have evolved according to the unperturbed system whose solution at time $t$ with initial condition $\eta$ is given by $\Th(t)\eta$. Applying the definition \eqref{eq:def:infgendef_illus} for this case we get the $\frac{d}{dt}|_{t=0}$ term in \eqref{eq:infgen_gennoise_L0_illus}. If there was noise alone we would get $\noisegent$ term in \eqref{eq:infgen_gennoise_L0_illus}. The rate of change of $\hprc^\eps$ in \eqref{eq:evolofhprceps_gennoise} is $\sigma\,\hprcdrift$ which explains the $\sigma \hprcdrift\, \frac{\partial f}{\partial h}$ term in \eqref{eq:infgen_gennoise_L1_illus}. The other term in \eqref{eq:infgen_gennoise_L1_illus} is due to the perturbation coefficient $\sigma F$ in \eqref{eq:detDDE_pert_rescale_gennoise}.

The problem of finding the infinitesimal generator $\igen_H$ of the process $\lim_{\eps \to 0}\hprc^{\eps}$ boils down to this (for details see the technique of martingale problem in chapter 5 of \cite{Ethier_Kurtz}): find an operator $\igen_H$ such that given any nice function $f_H$ of $h$ alone, there exists a function $f^\eps$ of $(\eta,\gnoise,h)$ such that $|f_H(h)-f^\eps(\eta,\gnoise,h)|$ and $|(\igen_Hf_H)(h)-(\igen^\eps f^\eps)(\eta,\gnoise,h)|$ are of order $\eps$.

Now we show how to find $\igen_H$. Formally, consider $f^\eps(\eta,\gnoise,h)\overset{\text{def}}=f_H(h)+\eps f_1(\eta,\gnoise,h)+\eps^2f_2(\eta,\gnoise,h)$ with $f_1$ and $f_2$ yet to be determined. Computing $\igen^\eps f^\eps$ we find
\begin{align}\label{eq:expandgenLepsFeps}
\igen^\eps f^\eps=\frac{1}{\eps^2}\igen_0f_H + \frac{1}{\eps}(\igen_0f_1+\igen_1f_0)+(\igen_0f_2+\igen_1f_1)+ O(\eps).
\end{align}
Note that $\igen_0f_H=0$ because $\igen_0$ involves differentials with respect to $(\eta,\gnoise)$ whereas $f_H$ is a constant as a function of $(\eta,\gnoise)$ (it is function only of $h$). Now, $f_1$ can be choosen so that $\igen_0f_1+\igen_1f_0=0$. It can be verified that $f_1$ is
\begin{align*}
f_1(\eta,\gnoise,h)=\int_0^\infty ds\left(\int_{\nstsp}\left(\nu(s,\gnoise,d\zeta)-\bar{\nu}(d\zeta)\right)\sigma(\zeta)\right)\hprcdrift(\Th(s)\eta)\frac{\partial f_H(h)}{\partial h}.
\end{align*}
We would not be able to select $f_2$ such that $\igen_0f_2+\igen_1f_1=0$. However $\igen_0f_2+(\igen_1f_1-\{\igen_1f_1\})=0$ can be solved where $\{\igen_1f_1\}$ is certain kind of average. With this choice of $f_2$, now \eqref{eq:expandgenLepsFeps} gives $|\igen^\eps f^\eps-\{\igen_1f_1\}|\sim O(\eps)$. Inspecting $\{\igen_1f_1\}$ gives $\igen_H$. Note that $\igen_1f_1$ equals
\begin{align*}
\int_0^\infty ds\left(\sigma(\gnoise)\int_{\nstsp}\left(\nu(s,\gnoise,d\zeta)-\bar{\nu}(d\zeta)\right)\sigma(\zeta)\right)\left((\Ind F(\eta).\nabla)\hprcdrift^{(s)}(\eta)\frac{\partial f_H(h)}{\partial h}+\hprcdrift(\eta)b^{(s)}(\eta)\frac{\partial^2f_H(h)}{\partial h^2}\right),
\end{align*}
where $\hprcdrift^{(s)}(\eta)\overset{\text{def}}=b(\Th(s)\eta)$. In the above expression (i) averaging the noise $\gnoise$ with respect to its invariant measure $\bar{\nu}$ and recalling the definition of autocorrelation in \eqref{eq:PREpap:autocorr}, (ii) realizing that $(\Ind F(\eta).\nabla)\hprcdrift^{(s)}(\eta)=(\Th(s)\Ind F(\eta).\nabla)\hprcdrift(\Th(s)\eta)$, and (iii) averaging the $\eta$ on trajectories of constant $h$, we get $\{\igen_1f_1\}$ as $b_H(h)\frac{\partial f_H(h)}{\partial h}+\frac12\sigma_H^2(h)\frac{\partial^2 f_H(h)}{\partial h^2}$ where $b_H$ and $\sigma_H$ are 
as stated in the theorem \ref{thm:gennoiseweakconvg}.


\section{\label{appsec:numericalscheme} Numerical scheme for simulations}
All simulations in this paper are done with Euler-Maruyama scheme. For example, \eqref{eq:examplesys_quadadded_forsimu} with $\gamma_c=0$ is simulated as follows. Select a time step $\Delta$. Let $N=r/\Delta$ where $r$ is the delay in the system. Specify initial conditions at the time points of the form $j\Delta$ for $j=-N,-N+1,\ldots,-2,-1,0$. Then, for $j\geq 0$,
$$x|_{(j+1)\Delta}=x|_{j\Delta}+\Delta \left(-\frac{\pi}{2}x +\eps \gamma_q x^2\right)\big|_{(j-N)\Delta} + \eps \sigma\sqrt{\Delta}\mathcal{N}_j,$$
where $\mathcal{N}_j$ is a standard normal random variable.

For \eqref{eq:detDDE_pert_gennoise_numsim} we first simulate the two-state markov chain and then use
$$x|_{(j+1)\Delta}=x|_{j\Delta}+\Delta \left(-\frac{\pi}{2} + \eps \sigma(\gnoise|_{j\Delta})\right) x\big|_{(j-N)\Delta}.$$

The following values of $\Delta$ are used: for section \ref{subsec:verifyNumSimMarkovChain} $\Delta=5\times 10^{-5}$, for section \ref{subsec:numverif} $\Delta=2\times 10^{-5}$, for section \ref{subsec:linearpertWhitenoise} $\Delta=10^{-5}$, for the stationary density in figure \ref{fig:Gaud_InvDens_with_Noise} $\Delta=5\times 10^{-6}$. 

When the delay system in section \ref{subsec:vdpolosc} is close to the stochastic bifurcation threshold, the probability density takes a long time to reach its steady state. Hence, using a small $\Delta$ was not practical. For example, generating figure \ref{fig:Gaud_InvDens_with_Noise} (which is not close to the stochastic threshold) for which the invariant density is reached by 4500 delay periods, took more than 24 hours on a Intel Xeon X5675 3.07GHz CPU with $\Delta=5\times 10^{-6}$ and $3200$ samples. Close to the bifurcation threshold, it takes much longer. Critical issue is not only the speed of CPU but also its memory. For simulating a delay system, the history of the process needs be stored in memory. Suppose one is simulating the stochastic delay system with 1000 samples; then, storing a $1000 \times (\frac{r}{5\times 10^{-6}})$ matrix where $r$ is delay in the system, requires a huge amount of memory. 

So, for computing stochastic bifurcation thresholds in section \ref{subsec:vdpolosc}, we used $\Delta=10^{-4}$ instead, at the expense of losing some accuracy. 

\vspace{12pt}

\end{document}